\documentclass[a4paper,10pt]{article}
 \pdfoutput=1
\usepackage[blocks]{authblk}
\usepackage{amsmath,amssymb,amsfonts,hyperref,comment,xcolor,url,dsfont}
\usepackage[T1]{fontenc}
\usepackage{fourier}
\usepackage{epsfig,graphicx,wrapfig,ccaption}
\usepackage{csquotes}

\usepackage{hyperref}
\usepackage{xcolor}

\definecolor{unbleu}{rgb}{0.03, 0.15, 0.4}
\definecolor{unvert}{rgb}{0.21, 0.37, 0.23}

\hypersetup{
pdfborder = {0 0 0},
colorlinks,
linkcolor=unbleu,
citecolor=unbleu,
urlcolor=orange
}
\makeatletter
\renewenvironment*{displayquote}
  {\begingroup\setlength{\leftmargini}{1cm}\csq@getcargs{\csq@bdquote{}{}}}
  {\csq@edquote\endgroup}
\makeatother

\begin{document}

\title{Computer experiments and visualization\\ in mathematics and physics\\ {\large A subjective short walk among some historical examples}}

\author[1]{J.-R. Chazottes
\thanks{Email: \texttt{jeanrene@cpht.polytechnique.fr}}}

\author[2]{M. Monticelli
\thanks{Email: \texttt{marc.monticelli@math.cnrs.fr}}}

\affil[1]{{\normalsize Centre de Physique Th\'eorique, CNRS, Institut Polytechnique de Paris, Palaiseau, France}}

\affil[2]{{\normalsize Laboratoire de Math\'ematiques J.A. Dieudonn\'e, Universit\'e C\^ote d'Azur, Nice, France}}

\date{\today}

\maketitle

\begin{abstract}
In this short essay, we show how computer experiments, and especially \emph{visualization}, allowed for the investigation and discovery of phenomena which would have 
passed unnoticed. We shall also highlight the importance of \emph{interactivity} between the computer and the user.
We do this by surveying several historical examples from mathematics and the physical sciences. Many pictures, and even hyperlinks to online interactive numerical 
experiments, are provided. Needless to say that we do not claim to be exhaustive, and that the chosen examples  reflect our taste as well as our limited knowledge. \newline
Hyperlinks to several online interactive digital experiments are provided.
\end{abstract}

\maketitle

\newpage 

\tableofcontents

%\newpage

%%%%%%% SECTION
\section*{Introduction}

The purpose of this article is to illustrate the impact of computer experiments on mathematics and the physical sciences from a historical perspective. 
We are especially interested in showing the role of visualization and interactivity. Let us start by giving the floor to John Hubbard, who will meet again in the section on
complex dynamics:
\begin{displayquote}
\textcolor{unbleu}{
``I mention computer graphics because faster and cheaper computers alone would not have had the same impact; without pictures, the information pouring out of 
mathematical computations would have remained hidden in a flood of numbers, difficult if not impossible to interpret. For people who doubt this, I have a story to relate. Lars 
Ahlfors, then in his seventies, told me in 1984 that in his youth, his adviser Lindelof had made him read the memoirs of Fatou and Julia, the prize essays from the Acad\'emie 
des Sciences in Paris. Ahlfors told me that they struck him at the time as ``the pits of complex analysis'': he only understood what Fatou and Julia had in mind when he saw the 
pictures Mandelbrot and I were producing. If Ahlfors, the creator of one of the main tools in the subject and the inspirer of Sullivan's no-wandering domains theorem, needed 
pictures to come to terms with the subject, what can one say of lesser mortals?''}
\end{displayquote}

In addition to their obvious brute-force capabilities of grinding out numerical solutions of differential equations, both John von Neumann and Stanislaw Ulam foresaw 
some of the more subtle applications of computers, as a flexible, interactive tool for the purpose of discoveries. 
%This interplay between computations and analysis, which Ulam called ``synergetics'', has indeed proved to be of great importance. 

After having evoked these two pioneers, we will outline the Fermi-Pasta-Ulam-Tsingou experiment which revealed a completely unexpected behaviour in 
a one-dimensional array of identical masses coupled to their nearest neighbors by slightly nonlinear springs (discretization of a vibrating string that included a small non-linear term).  \newline
Then we will see how Turing, motivated by his work on morphogenesis, suggested that numerical experiments should become a  genuine tool in scientific investigation. \newline
The next example will be taken from arithmetic geometry.
Birch and Swinnerton-Dyer used an early computer (which was the size of a large room and which had 20 kilobytes of memory!) to compute many examples of solutions to cubic equations in 
two variables modulo prime numbers. Graphing the output data in the right way led to new insights in the theory of elliptic curves which ultimately became the Birch and Swinnerton-Dyer 
conjecture, one of the Clay Millennium problems.\newline 
We will then discuss Lorenz's model which led him to observe the famous ``strange attractor'' bearing his name. \newline
Motivated by the Fermi-Pasta-Ulam-Tsingou experiment, Martin Kruskal and Norman Zabusky approximated the nonlinear spring-mass system by the 
Korteweg-de Vries equation, and numerically discovered ``solitons''. \newline
Our next example will be Michel H\'enon who was working in astrophysics. For instance, he was interested in the motion of a star around a galactic center 
for which he introduced a simplified model displaying a mixture of quasi-periodic and chaotic solutions. He also discovered a strange attractor bearing his name which is obtained by
iterating a simple two-dimensional nonlinear mapping.\newline
Then we will describe the period-doubling cascade and its universal features discovered independently by Pierre Coullet and Charles Tresser in Nice, and Mitchell Feigenbaum in Los Alamos. 
\newline
Our final example will be complex dynamics, with John Hubbard and Beno\^{\i}t Mandelbrot.

We do not claim to be exhaustive, and we are not historians. Among the noticeable people we could have mentioned there is Derrick H. Lehmer (1905-1991), a number 
theorist, who worked on the \textsf{ENIAC}. Fortunately, we can refer to the article of Liesbeth De Mol \cite{deMol} which 
discusses and contrasts the visions of Lehmer and von Neumann on the use of the computer and its impact on mathematics.

Notice that we deliberately avoid technicalities and refer, whenever possible, to some texts understandable by non-specialists. We also point to the original 
articles which are almost all freely available on Internet.

We did our best to find out who were the programmers behind the computer experiments, for instance Mary Tsingou for the FPUT model, Gary Deem who worked with
Kruskal and Zabusky, or Peter Moldave who wrote and run the computer experiments for Mandelbrot. We also did our best to indicate the computers used and show pictures of them
in the final section.

In this article we will not be concerned with exciting topics such as numerical analysis, computer-assisted proofs, proof assistants (formal proof management 
systems), or machine learning. 

Finally, let us mention that this article was originally supposed to be a translation of an article published in French in {\em La Gazette des Math\'ematiciens}, vol. 143 (2015).
In fact, the article was almost entirely rewritten and expanded, and many pictures were added.

{\em Acknowledgments.} We thank Alix Chazottes for a careful reading.

\newpage

%%%%%%% SECTION
\section{Von Neumann and Ulam, and ``synergesis''}\label{sec:vNU}

John von Neumann (1903-1957) and Stanislaw Ulam (1909-1984) seem to be the very first people to have understood the potential of computers
in mathematics and physics \cite{ulam}:
\begin{displayquote}
\textcolor{unbleu}{
``Almost immediately after the war Johny and I also began to discuss the possibilities of using computers
heuristically to try to obtain insights into questions of pure mathematics. By producing examples and by observing
the properties of special mathematical objects one could hope to obtain clues  as to the behavior of general statements
which have been tested on examples. In the following years in a number of published papers, 
I have suggested -- and in some cases solved -- a variety of problems in pure mathematics by such experimenting or
even merely ``observing'' \footnote{He of course means ``using a computer'' (A/N).}.''}
\end{displayquote}

Ulam proposed to use computers not as a ``mill'' producing numerical answers but as an experimental tool allowing
to study the solutions of the equations upon consideration. His vision had been influenced by physics, more
specifically by the numerical studies of simplified models of neutron diffusion (linked to the atomic bomb).\footnote{Let us mention in passing that the very first 
calculations for a chain reaction were made by Nicholas Metropolis in 1947 on \textsf{ENIAC} \cite{metropolis}.}

In a text of 1947 titled \textit{The Mathematician} \cite{vonNeumann}, von Neumann examined his work and asked himself whether 
mathematics was an empirical science. 
%His attempt at answering such a question is by a comparison with theoretical physics.
For him, there is an empirical basis to mathematics which was overshadowed by their axiomatical development but,
\textcolor{unbleu}{``when it shows signs of becoming baroque, then the danger signal is up''} and \textcolor{unbleu}{``when this stage is reached, the only remedy seems to me the 
rejuvenating return to the source: the reinjection of more or less directly empirical ideas.''}

\begin{figure}[htb!]
\centering
 \includegraphics[width=0.49\textwidth]{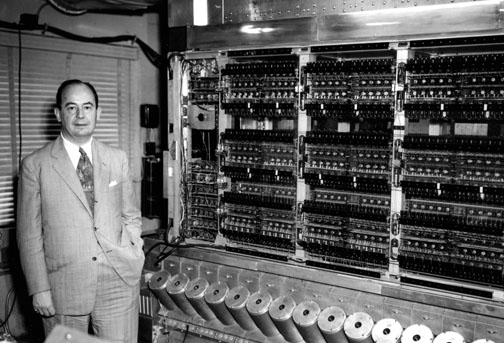}
\includegraphics[width=0.42\textwidth]{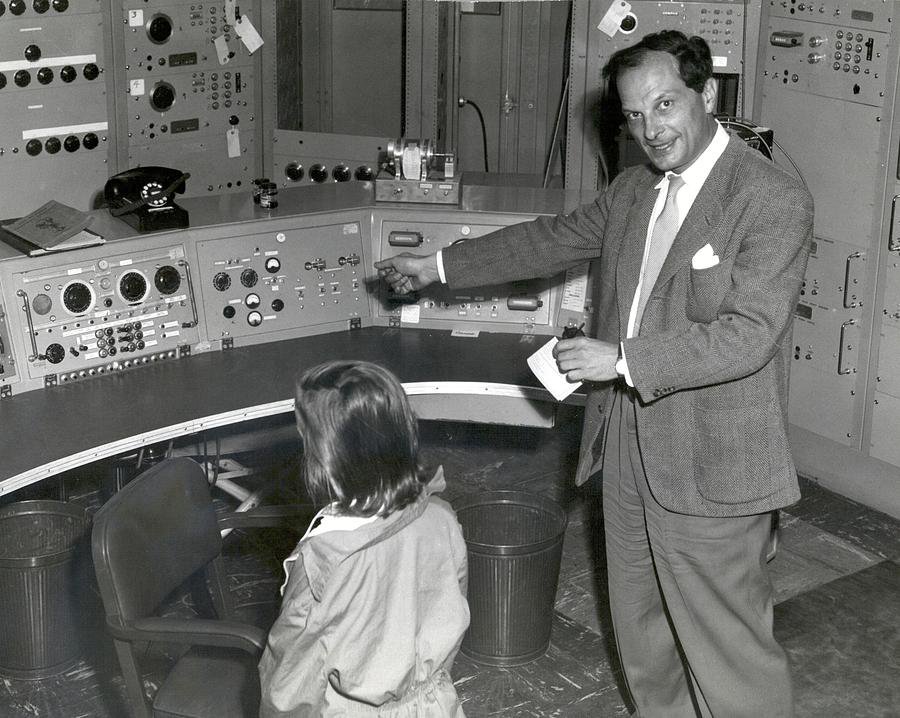}
\legend{\scriptsize{von Neumann and the IAS machine (1952) and S. Ulam at MANIAC panel with his daughter Claire (1955)}}
%\label{} % a mettre apres caption
\end{figure}

\begin{figure}[htb!]
\centering
\includegraphics[scale=.14]{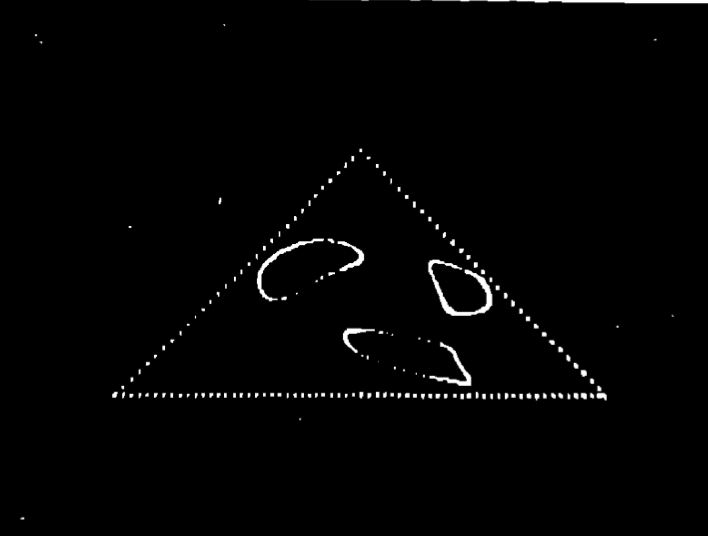}\hspace{0.0cm}
\includegraphics[scale=.15]{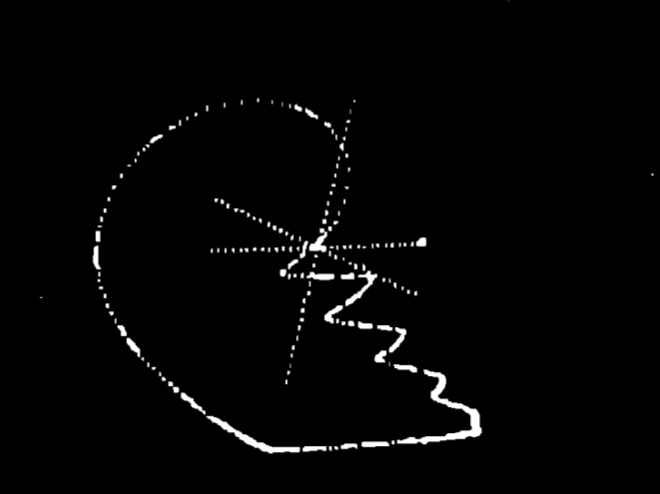}\\
\includegraphics[scale=.1375]{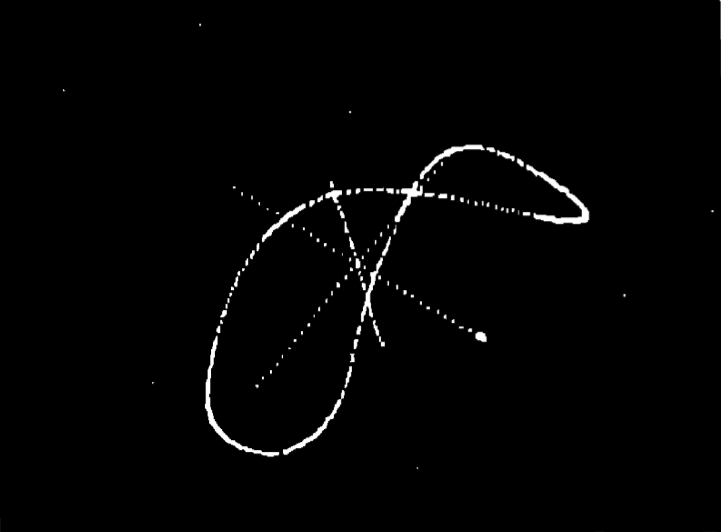}\hspace{0.01cm}
\includegraphics[scale=.149]{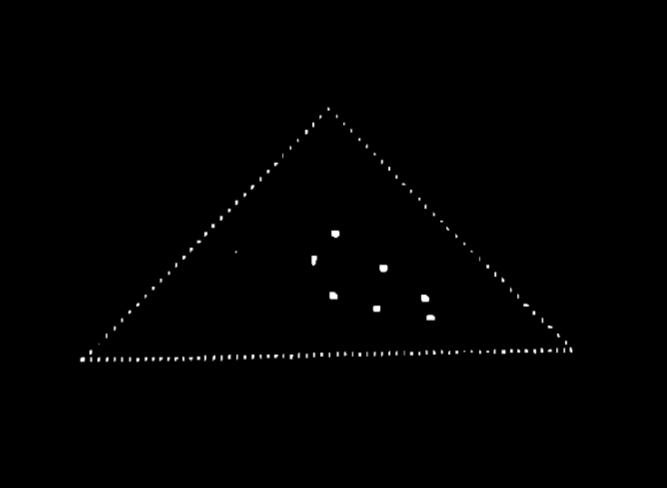}
\legend{Oscilloscope (!) plots from an article of Stein and Ulam \cite{stein-ulam}.}
\label{fig-ulam-stein} % a mettre apres caption
\end{figure}

\begin{figure}
\centering
\includegraphics[scale=.2]{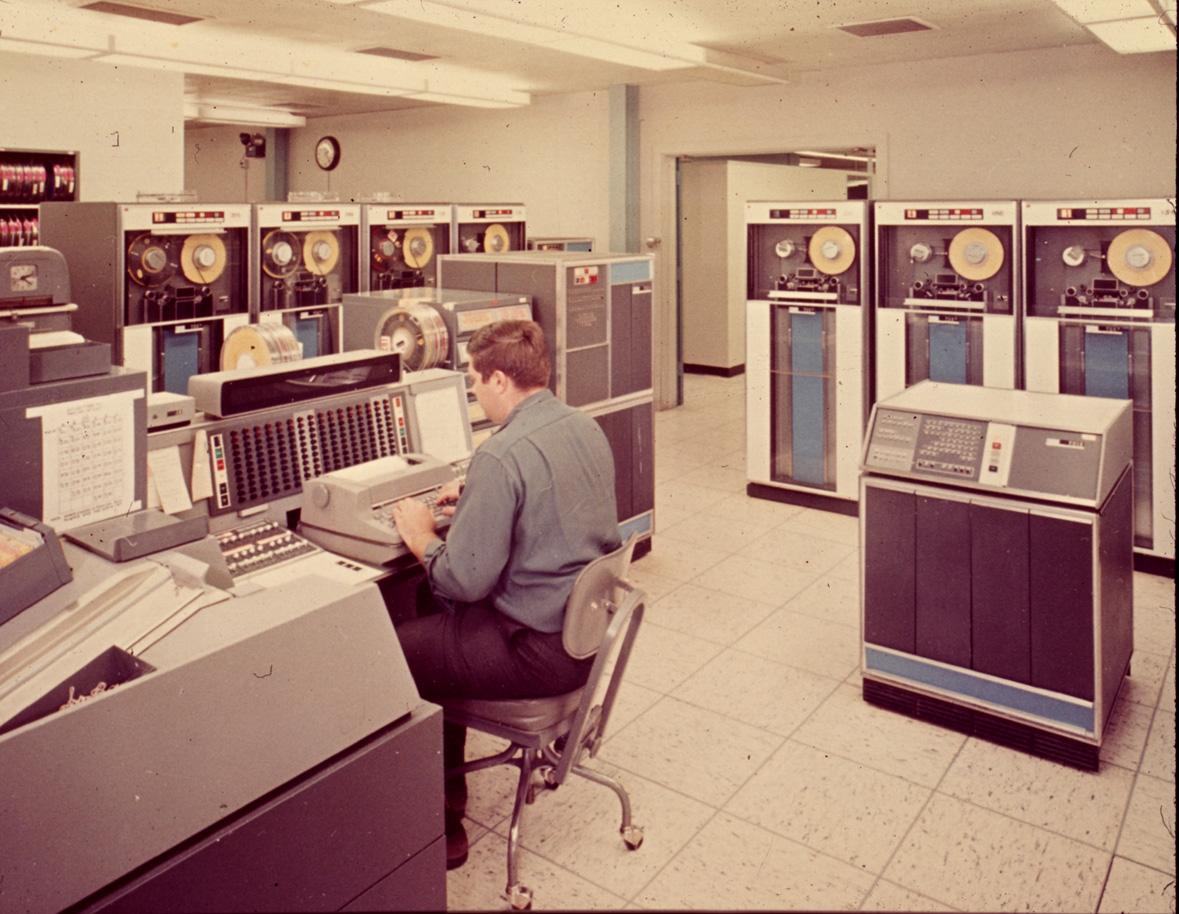}
\legend{The IBM 7030, also known as Stretch, was IBM's first transistorized supercomputer. It appeared in 1961.}
\label{fig-IBM7030} % a mettre apres caption
\end{figure}

In this spirit, Ulam aimed at developing a practice of numerical experiments, starting with combinatorics and number theory. Then he turned with his 
collaborators to the exploration of the behavior of the iterations of nonlinear transformations which are the counterpart in discrete time to ordinary 
differential equations \cite{stein-ulam}.\footnote{We cannot resist quoting the MathSciNet review for this article: ``   The paper considers mainly the three-dimensional iteration transformations 
$x_i=P_i(x_1,x_2,x_3)$, $i=1,2,3$ where the $P_i$ are cubics in $x_1,x_2,x_3$. [...] Many photographs of cathode ray tube displays are given, a fondness for citing large numbers of iterations 
and machine time used is revealed, and a crude classification of the limited results is offered, but there appears to be no firm new results of general mathematical interest. It 
is an attractive idea that computers can produce many specific cases of mathematical situations from which we can obtain insight, but usually (though not always) it is 
the insight, not the specific cases, that has been the criterion 
for public presentation. One can only wonder what will happen to mathematics if we allow the undigested outputs of computers to fill our literature. The present paper 
shows only slight traces of any digestion of the computer output.''}
They use a device (an oscilloscope!) connected to the computer which allows the visualization of orbits -- a new way to study nonlinear iterations that seems so obvious 
nowadays. An example of this approach is found in the paper \cite{stein-ulam} from which we extracted a few plots; See Fig. \ref{fig-ulam-stein}. They used an
\textsf{IBM 7030 ``STRETCH''} available at the Los Alamos Scientific Laboratory; see Fig. \ref{fig-IBM7030}. (All numbered figures are gathered in Section \ref{sec:more-pics}.)

In his book \cite{ulam-collection}, published in 1960, Ulam clearly outlines what we nowadays call ``interactivity'' (between the user and the computer), which he calls ``synergesis''. In 
the chapter titled \textit{Computing Machine as a Heuristic  Aid}, he writes: 
\begin{displayquote}
\textcolor{unbleu}{`` Instead of using the machine as a robot or, as it were, as a player piano whose tunes are written 
in advance, the machine kept in constant communication with an intelligent operator who changes even the logical nature of the problem during the course of a computation, at will, after 
evaluating the results which the machine provides. Of course such possibilities exist already, but to a \emph{very limited}\,\footnote{Ulam's emphasis} extent. [...]
Obviously for activites of this sort, a rapid access to the machine is necessary [...]. Also, the machine has to provide a quick illustration and display of the computed quantities and figures.''
}
\end{displayquote}

%%%%%%% SECTION
\section{Fermi, Pasta, Ulam and Tsingou: The FPUT experiments on \textsf{MANIAC}}\label{sec:FPUT}

Let's Ulam open this section :
\begin{displayquote}
\textcolor{unbleu}{
``Computers were brand-new; in fact the Los Alamos Maniac was barely finished. The Princeton von Neumann machine had met with technical and engineering difficulties that had
prolonged its perfection. [...] As soon as the machines were finished, Fermi, with his great common sense and intuition, recognized immediately their importance for the study of
problems in theoretical physics, astrophysics, and classical phy\-sics. We discussed this at length and decided to attempt to formulate a problem simple to state, but such that a solution
would require lengthy computation which could not be done with pencil and paper or with existing mechanical computers. [...] we found a typical one requiring long-range prediction and
long-time behavior of a dynamical system. It was the consideration of an elastic string with two fixed ends, subject not only to the usual elastic force of strain proportional to strain, but having, in
addition, a physically correct small non-linear term. The question was to find out how this non-nonlinearity after very many periods of vibrations would gradually alter the well-known
periodic behavior of back and forth oscillation in one mode; how other modes of the string would become more important; and how, we thought, the entire motion would ultimately thermalize [...] John Pasta, 
a recently arrived physicist, assisted us in the task of flow diagramming, programming, and running the problem on the Maniac. Fermi had decided to try to learn how to code the machine by himself. [...]
Our problem turned out to have been felicitously chosen. The results were entirely different qualitatively from what even Fermi, with his great knowledge of wave motions, had expected. 
} (Quoted from \cite[Chapter 12]{ulam-autobio}.)
\end{displayquote}

\begin{wrapfigure}{r}{0.48\textwidth}
\vspace{-10pt}
  \begin{center}
    \includegraphics[width=0.2\textwidth]{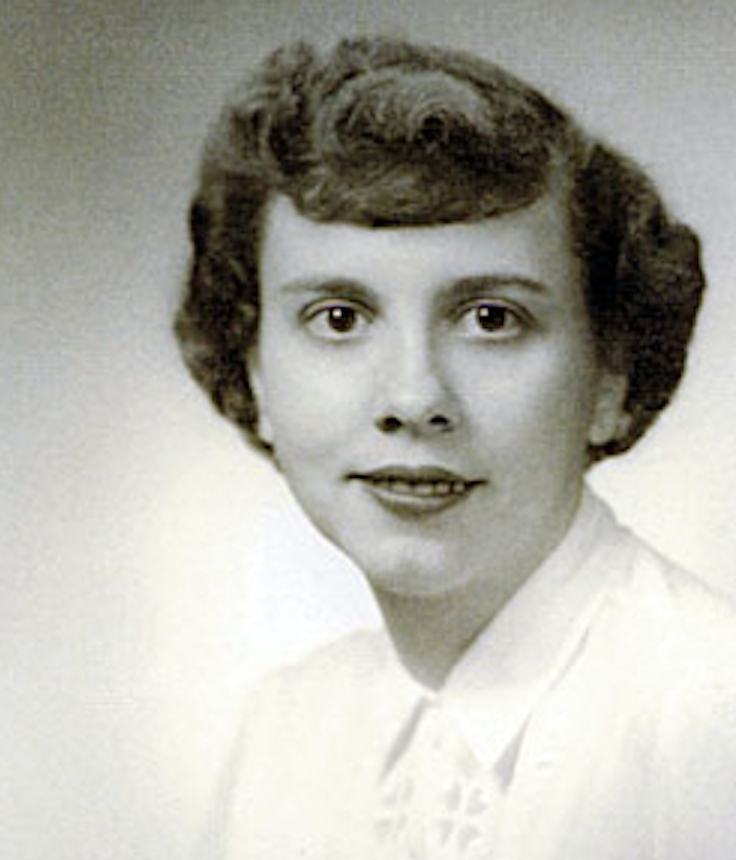}
    \includegraphics[width=0.199\textwidth]{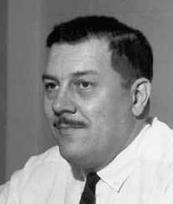}\\
    \includegraphics[width=0.405\textwidth]{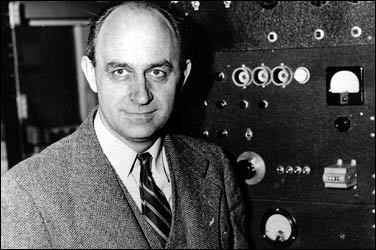}
  \end{center}
  \vspace{-15pt}
  %\caption{\scriptsize{}}
  \legend{{\footnotesize Mary Tsingou, John Pasta and Enrico Fermi}}
  \vspace{-10pt}
\end{wrapfigure}
More precisely, the model, obtained as the discretization of a partial 
differential equation model of a string, consists in a one-dimensional array of identical masses coupled to their nearest neighbors by springs, with fixed end-points. 
When a mass is moved away from its position at rest, it undergoes a force pulling it back to equilibrium.\footnote{The equations are the following:
\[
\ddot{u}_j=u_{j+1}-2u_j+u_{j-1} +\alpha\Big[\big(u_{j+1}-u_j\big)^2-\big(u_j-u_{j-1}\big)^2\Big]
\]
where $u_j$ is the relative displacement with respect to the equilibrium position of the $j$-th mass, $u_0=u_N=0$ (ends of the chain assumed to be fixed), 
and $\alpha$ is the `small' nonlinearity parameter ($\alpha=0$ gives back the harmonic oscillator).}
They considered 16, 32 and 64 masses. To their great surprise, the system did not exhibit thermalization but instead a complex quasi-periodic behavior. This was in contradiction 
with the so-called ``ergodic hypothesis'' which was assumed to hold true in this case.
There are rumors throughout the literature that the \textsc{MANIAC I} was ``accidently'' left on one night and that the scientists came back after realizing their mistake to find astonishing results. According to Menzel\footnote{Telephone interview realized by L. Freeman in her Bachelor thesis \cite{freeman}.}, the computers were run at night because the numeric iteration computations for long time scales were very slow and the computer was being used for weapons design during the day.

A technical report appeared in 1955 \cite{FPU}\footnote{Fermi died in 1954 and this report was published some ten years later as part of Fermi's collected works, and
were never published in a journal.} but Mary Tsingou (born in 1928) who wrote the algorithm and programmed \textsf{MANIAC I} is not on the list of authors.
In a footnote it is written: ``We thank Miss Mary Tsingou for efficient coding of the problems and for running the computations on the Los Alamos MANIAC machine.''
Notice that Tsingou is not mentioned in Ulam's quote above (and in fact nowhere in his book). Nowadays, this model is named the FPUT model instead of the FPU model. We refer to \cite{DPR} for more informations.

We cannot resist quoting what Pasta says about \textsf{MANIAC I} at a conference held in 1977: 
\begin{displayquote}
\textcolor{unbleu}{``The program was of course punched on cards. A DO loop was executed by the operator feeding in the deck of cards over and over again until
the loop was completed!''.}
\end{displayquote}
We notice that there were people who believed that thermal relaxation times for the FPUT system were 
much too long to be observed during the short integration runs made by Tsingou, in other words they believed that the observed behavior was a transient effect. But Tsingou and Tuck 
\cite{tuck-menzel} made this conjecture very unlikely by  extensive numerical experiments in 1972.\footnote{By that time she was going by her married name Menzel.}

%%%%%%% SECTION
\section{Turing, morphogenesis, and computers}\label{sec:T}

Alan Turing (1912-1954) is best known for the creation of the ``Turing machine'' which he used to solve, in a masterful way, Hilbert's decision problem, and also for helping 
crack the code on intercepted Nazi messages helping the Allies win many major engagements during World War 2. 
Turing can also be considered a forefather of artificial intelligence with his seminal paper ``Computing Machinery and Intelligence'' published in 1950. 
In this section we give a very brief account on his pioneering work on morphogenesis. For an account of most Turing works, we refer to \cite{turing-guide}.
\begin{wrapfigure}{r}{0.48\textwidth}
\vspace{-10pt}
  \begin{center}
    \includegraphics[width=0.28\textwidth]{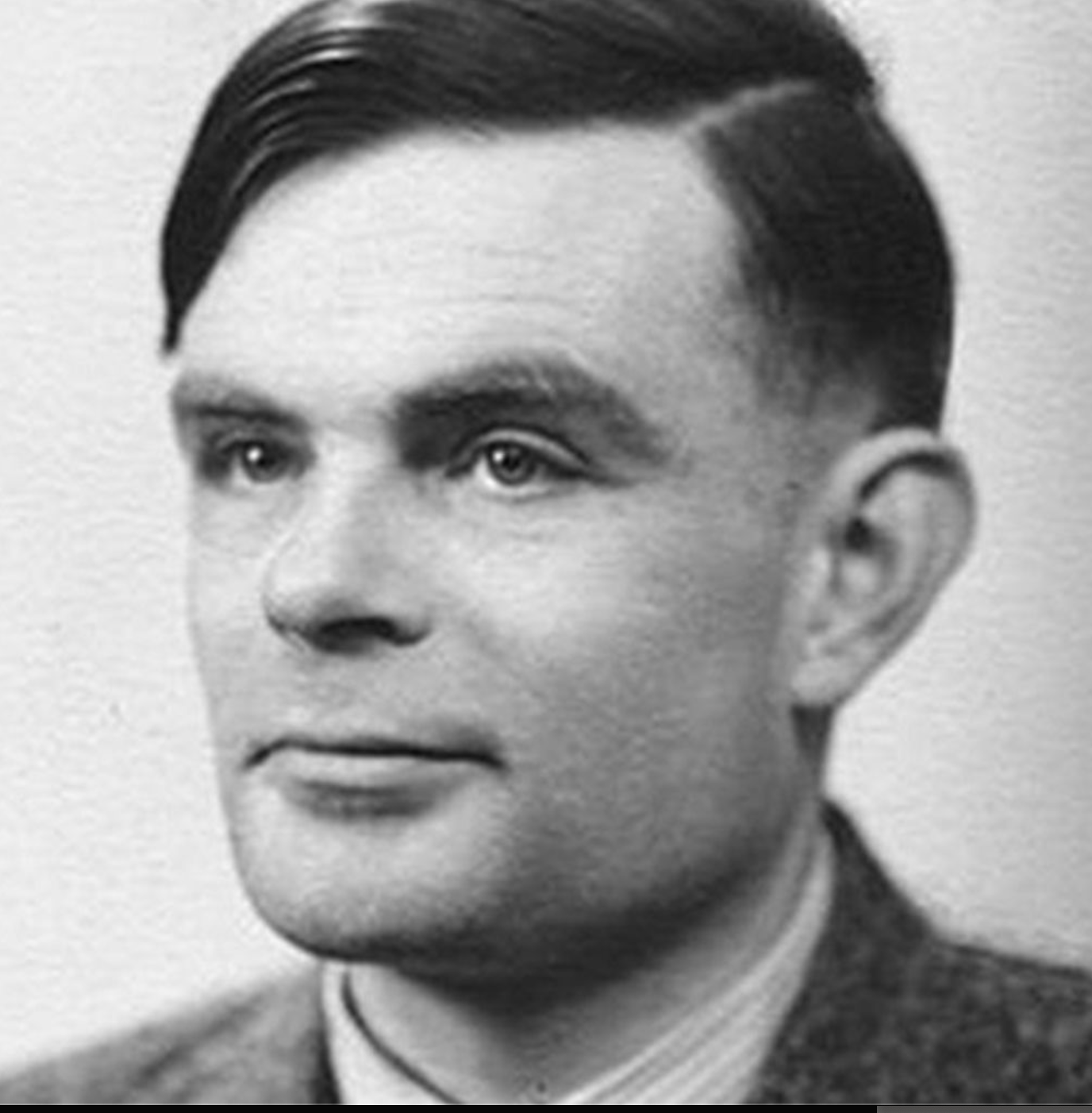}
  \end{center}
  \vspace{-15pt}
  %\caption{\scriptsize{}}
  \legend{{\footnotesize Alan Turing}}
  \vspace{-16pt}
\end{wrapfigure}

One of the questions that the biologist D'Arcy Thompson had asked himself was the emergence of similar forms for
unrelated organisms, making them inexplicable by purely genetic factors.
Turing postulated that there must be some general underlying process obeying physicochemical laws. 
He worked on setting up a mathematical model whose purpose was to account for the ``morphogenesis'', that is, the transition from an initial symmetric equilibrium state to a 
new equilibrium state breaking the initial symmetry and giving a form. This transition was modeled as a ``reaction-diffusion'' process within the chemical components 
of the system.

In 1952, Turing published a seminal paper titled ``The chemical basis of morphogenesis'' \cite{turing}\footnote{This paper contains only six references, one of them being the 
famous book of D'Arcy Thomson, ``On Growth and Form''.} in which he described his model and discussed two examples:
\begin{itemize}
\item
the development of spots like the ones appearing on the pelage of certain animals like leopards;
\item
the freshwater polyp Hydra. Initially its tube-shaped body is symmetric until a head with between five and ten tentacles appears at one end of the 
body.
\end{itemize}
The general model he proposed describes the interaction between two chemicals he calls ``morphogens''.
One is an ``activator'', which is autocatalytic and so introduces positive feedback. The other is an ``inhibitor'', which suppresses the 
autocatalysis of the activator. Crucially, they must have different rates of diffusion, the inhibitor being faster. In effect, this means that the activator's 
self-amplification is corralled into local patches, while the inhibitor prevents another such patch from growing too close by.

\begin{figure}[htb!]
\centering
\includegraphics[angle=0.8,width=0.4\textwidth]{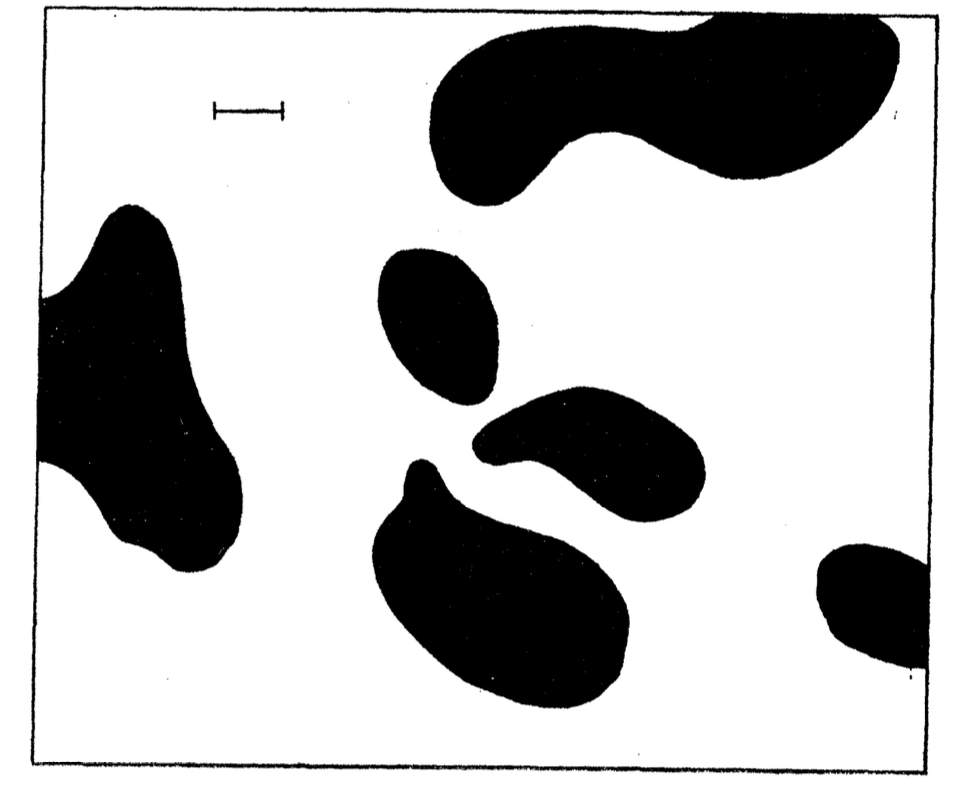}\hspace{0.5cm}
\includegraphics[angle=0.8,width=0.4\textwidth]{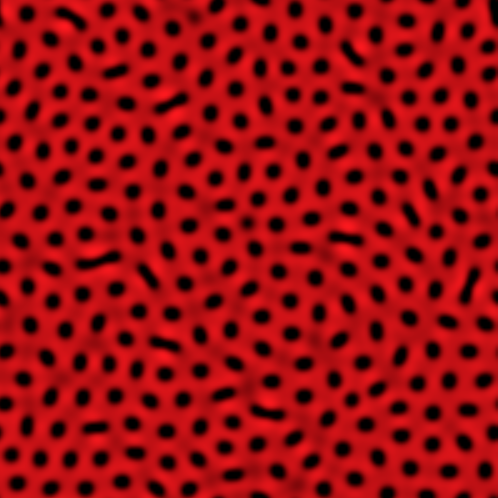}
%\caption{\scriptsize{}}
\legend{{\footnotesize Left: Dappling produced by the reaction-diffusion equations; this picture is taken from Turing's article. It was
hand-drawn by him on the basis of numerical calculations. Right: A Turing pattern obtained by the authors. The equations are
$\partial u/\partial t=u(v-1)+A \nabla^2\! u, \partial v/\partial t=16-uv+B \nabla^2\! u$ where $u=u(x,y,t)$ is the concentration of the activator at
point $(x,y)$ at time $t$, $v=v(x,y,t)$ that of the inhibitor, and $A,B$ are the diffusion coefficients.
\href{https://experiences.mathemarium.fr/Turing-patterns.html}{\textsc{click here}} for an online interactive digital experiment to otbain the above pattern and many others.}}
\end{figure}

It turns out that Turing made most of his calculations by hand but he shows a numerical example realized on \textsf{Manchester Mark I}.
At the end of his article, Turing suggests that numerical experiments should become a genuine tool in scientific investigation.
He writes: 
\begin{displayquote}
\textcolor{unbleu}{
``The difficulties are, however, such that one cannot hope to have any very embracing \emph{theory}  of such processes, beyond the statement
of the equations. It might be possible, however, to treat a few particular cases in detail with the aid of a digital computer. This method has the 
advantage it is not so necessary  to make simplifying assumptions as it is when doing a more theoretical type of analysis.''}
\end{displayquote}

Less known is his work on the problem of phyllotaxis (that is, the arrangement of leaves on a plant stem), which was never published at his time but is included in 
Turing's collected works \cite{turing-collected-works}. In that work, he indeed used computer simulations. 

Let us mention that it is only in the 1990s that  ``Turing patterns'' were obtained for the first time in a chemistry
experiment \cite{turing-exp}.

We close this section by mentioning that, in March 1946, Turing presented the world's first complete design for a stored-program electronic computer, ACE. Although Turing had seen the draft report on EDVAC by von Neumann, the ACE design was very different and included detailed circuit diagrams as well as software examples, and a precise budget estimate. Unfortunately, this projectdid did not go as planned.

%%%%%%% SECTION
\section{Birch and Swinnerton-Dyer: the rational points of elliptic curves}\label{sec:BSD}

%the rank of the elliptic curve is the number of independent rational points of infinite order

\begin{wrapfigure}{r}{0.48\textwidth}
\vspace{-10pt}
  \begin{center}
    \includegraphics[width=0.4\textwidth]{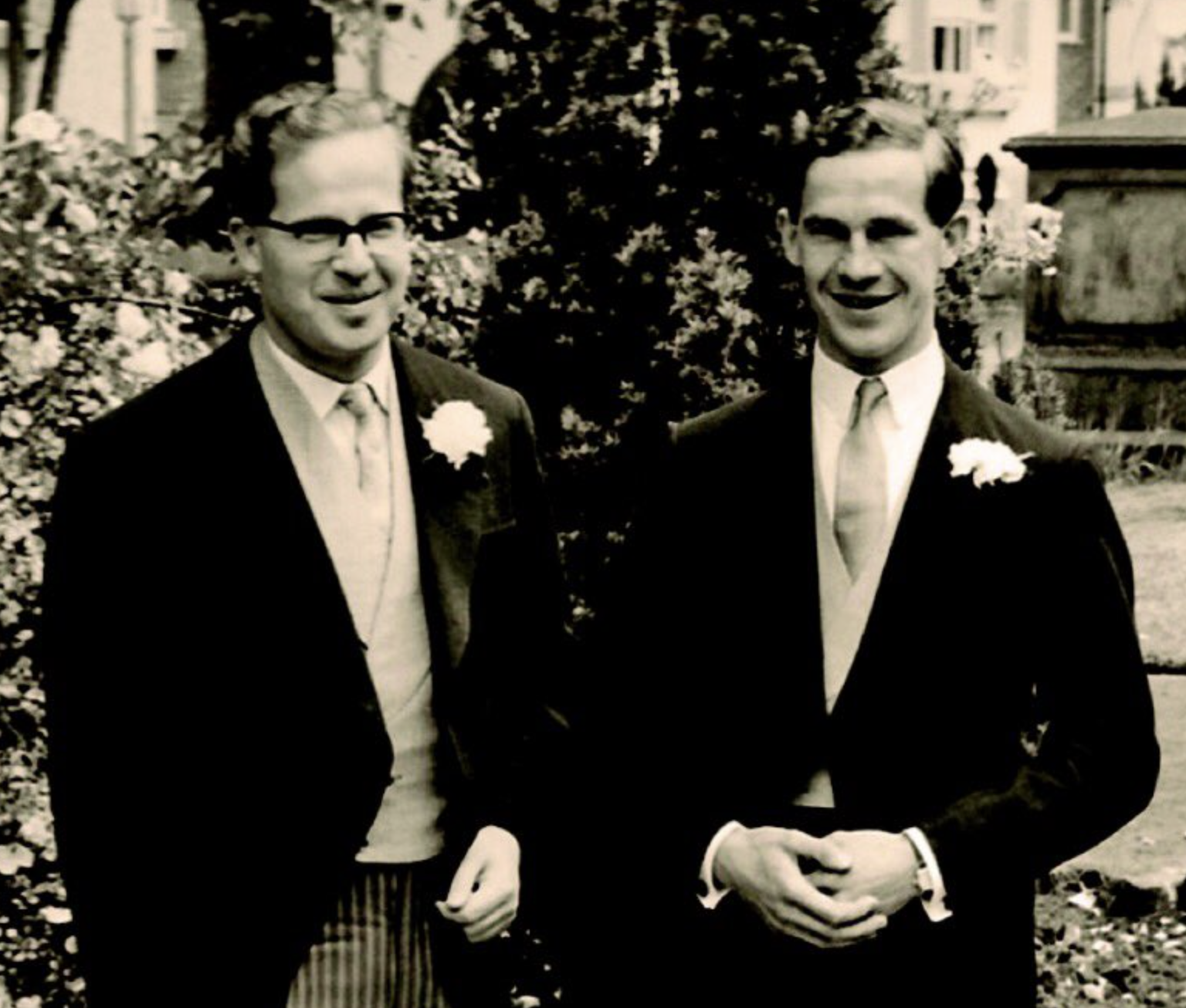}
  \end{center}
  \vspace{-15pt}
  %\caption{\scriptsize{}}
  \legend{{\footnotesize Peter Swinnerton-Dyer and Bryan Birch}}
  \vspace{-10pt}
\end{wrapfigure}

In this section, we will describe very roughly the Birch and Swinnerton-Dyer conjecture\footnote{As of 2021, only special cases of the conjecture have 
been proven.} which is related to the number of points with rational coordinates on ``elliptic curves''. This amounts to solving a special class of Diophantine equations, that is, polynomial 
equations whose coefficients are integers or rational numbers, and looking for their integer or rational solutions. 
The conjecture was developed at the beginning of the 1960s by the mathematicians Bryan 
Birch (born in 1931) and Peter Swinnerton-Dyer (1927-2018) with the help of \textsf{EDSAC} (an \textsf{ENIAC} descendant), at University of 
Cambridge Computer Laboratory. An elliptic curve $E$ is the set of solutions to an equation of the form $y^2=x^3+ax+b$, where $a$ and $b$ are
rational numbers such that $4a^3+27b^2\neq 0$.\footnote{This condition is equivalent to the curve being smooth, for instance this rules out cusps.}

\begin{wrapfigure}{r}{0.5\textwidth}
\vspace{-20pt}
  \begin{center}
    \includegraphics[width=0.4\textwidth]{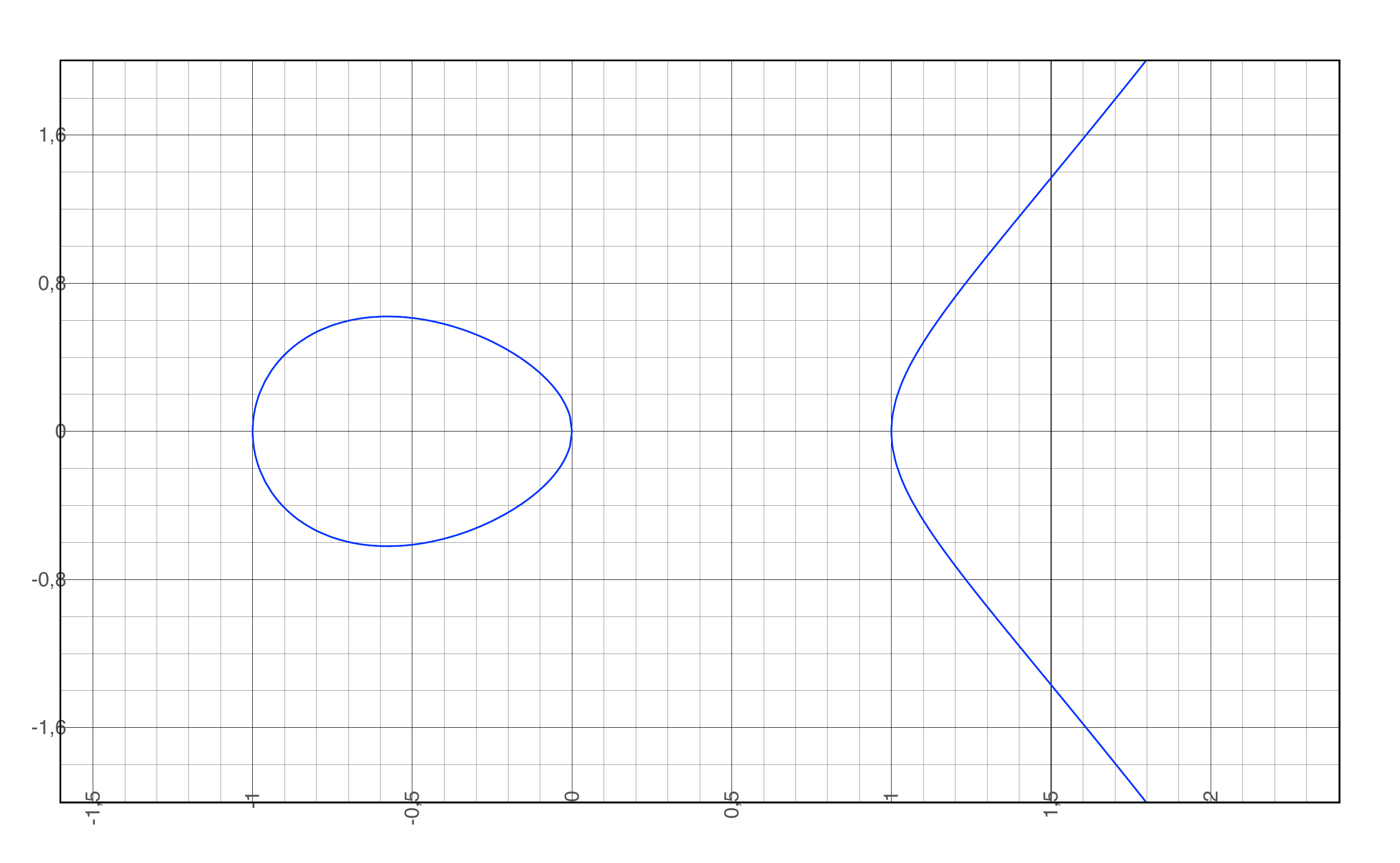}
  \end{center}
  \vspace{-20pt}
  %\caption{\scriptsize{}}
  \legend{{\footnotesize Elliptic curve defined by $y^2=x^3-x$.}}
  \vspace{-10pt}
\end{wrapfigure}

On the opposite figure, we show an example with $a=-1$ and $b=0$. The point $(2,\sqrt{6})$ lies on $E$ but it is not an integer-valued point. By setting 
$y=0$, we obtain the integer-valued points $(0,0)$, $(1,0)$ and $(-1,0)$. Fermat proved by his method of infinite descent that these three points are 
the only rational points on $E$.

A striking feature of rational solutions to elliptic curve equations is that these solutions form a group, thanks to a geometric construction we do not 
describe here, see {\em e.g.} \cite{rubin-silverberg}. The resulting structure is called the Mordell-Weil group of the curve, which is finite, and of course 
number theorists want to calculate it. That involves  finding a system of generators: rational solutions from which all others can be deduced by 
repeatedly using the group operation. At least, we would like to know how big this group is. By Mordell's theorem, this group can be broken into
two disjoint pieces that are themselves groups: An infinite part, which can be broken into a finite number of disjoint copies of the group of
integers $\mathds{Z}$, and a finite part, which can be broken into a finite number of disjoint copies of the group $\mathds{Z}/m\mathds{Z}$ of
integers modulo $m$, for certain integers $m$. The finite-part piece is rather well understood. The difficult part, which is the subject of the 
Birch and Swinnerton-Dyer conjecture, is to find the number of copies of $\mathds{Z}$ that appear in the infinite-part piece. This number is
called the ``rank'' of the elliptic curve, and it is denoted by $r$.\footnote{To get a group, one has to add a point ``at infinity'' to the rational points of the 
given elliptic curve $E$, because the ``right'' way to to view an elliptic curve is to view it as a curve in projective space $\mathds{P}^2$. Then 
this group can be written as the direct sum of $\mathds{Z}^r$ and a ``torsion subgroup'' which is characterized by a theorem of Mazur; see \cite{rubin-silverberg}.}

Let us come back to Birch and Swinnerton-Dyer. Their idea was to numerically look for the points on an elliptic curve modulo a given prime $p$, as 
there are only a finite number of possibilities to check. Of course, for large primes it is computationally intensive. They computed the number of points
modulo $p$ (denoted by $N_p$) on the elliptic curves $y^2=x^3-dx$, for five values of $d$ for which the rank was known. 
Birch noticed that Swinnerton-Dyer's computer experiments produce an interesting pattern if you divide $N_p$ by the prime concerned. Then multiply 
all of these fractions together, for all primes less than or equal to a given one, and plot the results against successive primes on logarithmic graph 
paper. The data all seem to lie close to a straight line, whose slope is the rank of the elliptic curve.
They obtained the graphical plots shown in the above figure.

\begin{figure}
\begin{center}
\includegraphics[width=0.6\textwidth]{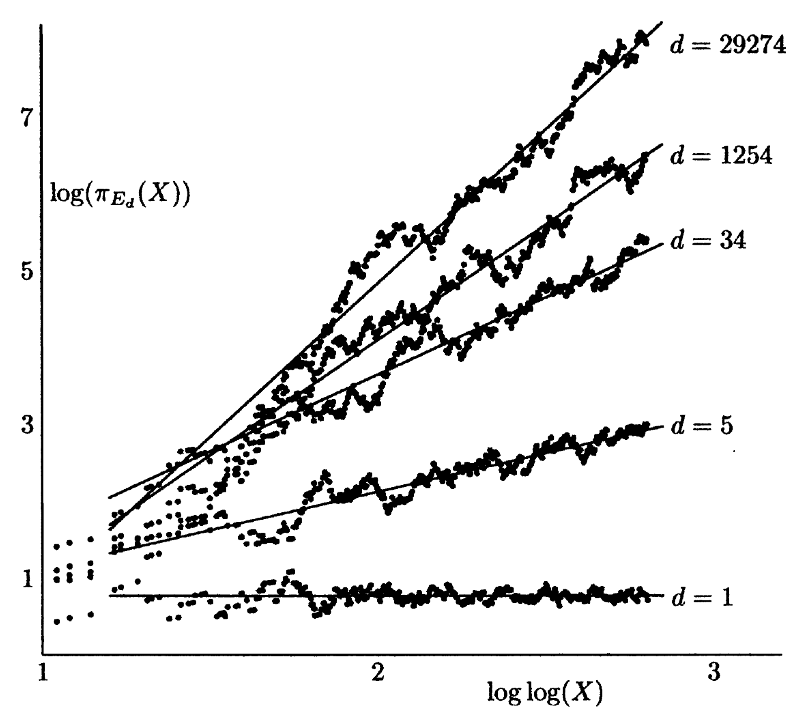}
\end{center}
\legend{{\footnotesize From \cite{BSD}. One has $r=0$ for $d=1$, $r=1$ for $d=5$,
$r=2$ for $d=34$, $r=3$ for $d=1254$, and $r=4$ for $29274$.}}
\end{figure}

They conjectured that
\[
\pi_{E_d}(X):=\prod_{p\leq X}\frac{N_p}{p} \sim C\, (\log X)^r,\quad X\to\infty,
\]
where $C$ is a constant. In fact, the function $\pi_E$ is difficult to work with, so they stated a related conjecture involving the so-called
$L$-function of $E$ in place of $\pi_E$. We will not give any detail here and again refer to \cite{rubin-silverberg}.

Let us finish this section by quoting Ian Stewart \cite{stewart} :
\begin{displayquote}
\textcolor{unbleu}{
``In the 1960s, when computers were just coming into being, the University of Cambridge had one of the earlier ones, called \textsf{EDSAC}. 
Which stands for electronic delay storage automatic calculator, and shows how proud its inventors were of its memory system, which sent sound 
waves along tubes of mercury and redirected them back to the beginning again. It 
was the size of a large truck, and I vividly remember being shown round it in 1963. Its circuits were based on thousands of valves - vacuum tubes. 
There were vast racks of the things along the walls, replacements to be inserted when a tube in the machine itself blew up. Which was fairly often.''
}
\end{displayquote}

%%%%%%% SECTION
\section{Lorenz: from meteorology to strange attractors}\label{sec:Lorenz}

At MIT, in 1953, Edward Lorenz was put in charge of a project devoted to statistical forecasting, exploring how the newly available digital computer could be put to use. 
A key issue was to know how well such numerical tools, based on linear statistical models, could predict complex weather patterns. Lorenz was
skeptical of the appropriateness of such tools as most atmospheric phenomena are non-linear. He sought a set of simple nonlinear differential equations that would 
mimic meteorological variations, thus providing a test example for the linear statistical approach. It was necessary to find a system with non-periodic
behavior (think of the large-scale turbulent eddies such as cyclones and anticyclones). After several failed attempts and some success with more complicated models,
he arrived in 1961 at the following system of ordinary differential equations:
\[
\begin{cases}
\dot{x} =\sigma (y-x)\\
\dot{y}= rx-y-xz\\
\dot{z}=xy-bz
\end{cases}
\]
where $\dot{x}=\mathrm{d} x/\mathrm{d} t$, etc, and where $\sigma$ is the Prandtl number (ratio of fluid viscosity to thermal conductivity), $r$ represents a temperature 
difference driving the system, and $b$ is a geometrical factor. What matters for us is that Lorenz numerically solved them and made a 
breakthrough: he discovered what was then  called ``deterministic chaos''. (In his paper, he used the values $\sigma=10$, $r=28$, $b=8/3$.)
The crucial insight of Lorenz was to observe that although orbits do depend on initial conditions, they accumulate on a kind of ``surface'' with ``figure-eight'' shape
which is insensitive to initial conditions. 

Lorenz made a rough sketch of this object and understood that two close initial conditions had the property of rapidly converging toward this ``surface'' and of travelling together for a while 
after which they start separating, at seemingly random intervals -- one staying in a ``wing'' while the other goes to the other one -- before they come close to 
each other back again, and so forth. 
Lorenz had just discovered what David Ruelle and Floris Takens later called a ``strange attractor'' \cite{ruelle-takens}. Strange because
this is not a fixed point or a closed orbit (limit cycle), and not generally a manifold. It turns out to be locally the product of a Cantor set and
a piece of two-dimensional manifold. 

\begin{figure}[htb!]
\centering
\includegraphics[width=0.24\textwidth]{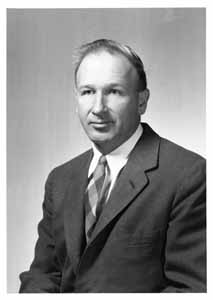} \hspace{.1cm}\includegraphics[width=0.3\textwidth]{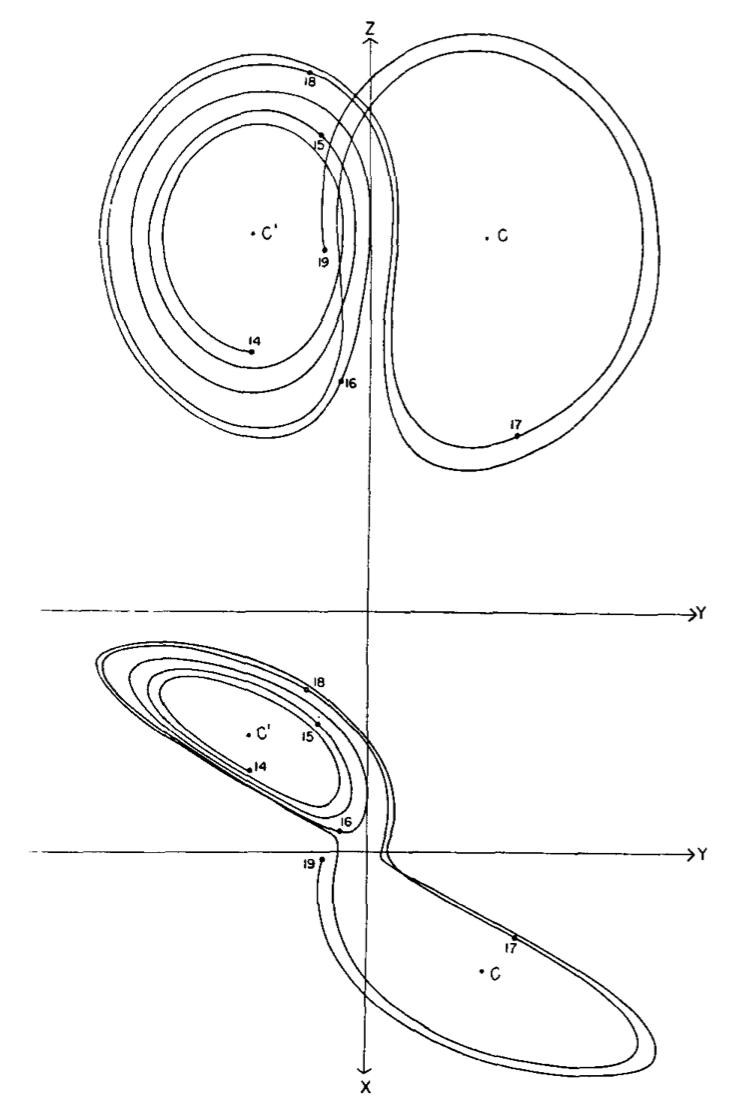}  \hspace{.1cm}\includegraphics[width=0.42\textwidth]{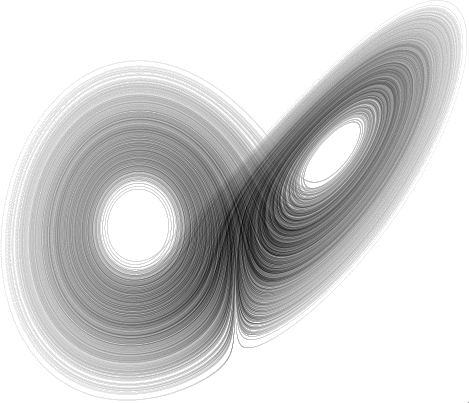}  
\legend{{\footnotesize Lorenz (in 1956), a figure extracted from his paper, and his attractor (simulation by the authors).\newline
\href{https://experiences.mathemarium.fr/Lorenz-attractor.html}{\textsc{click here}} for online interactive digital experiments.}}
\end{figure}

Let us emphasize Lorenz's lucidity and stroke of genius in not attributing what he observed to an artifact of the
computer or an effect specific to his model. Instead he understood that there was an underlying general phenomenon.
Indeed, one should realize that Lorenz used a \textsf{Royal McBee LGP-30}. 
Margaret Hamilton (born in 1936) and then Ellen Fetter (born in 1940, also known as Ellen Gille after she married) were Lorenz's programmers. 
As all computers of that time, it was slow, noisy and not as reliable as today's 
computers.

\begin{figure}[htb!]
\centering
   \includegraphics[width=0.25\textwidth]{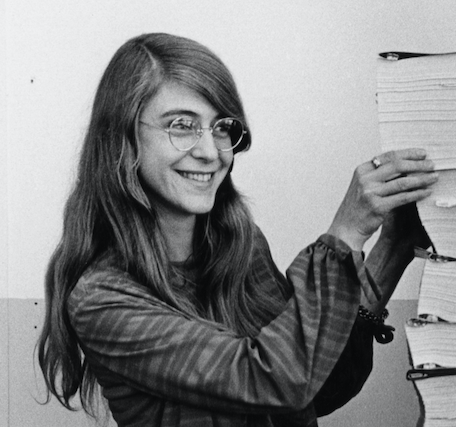} \hspace{1cm}\includegraphics[width=0.175\textwidth]{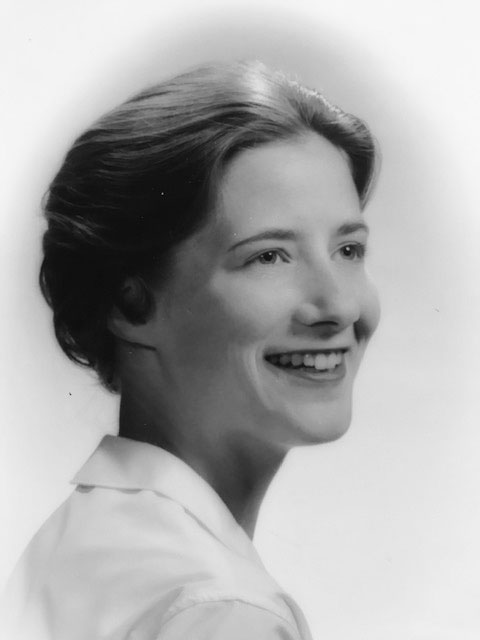}  
  \legend{{\footnotesize M. Hamilton and E. Fetter in the 1960s.}}
\end{figure}

%\textbf{This ``desk'' computer -- it was the size of a desk -- weighed some 350kg and sounded like a passing propeller plane. It was so loud that it even got its own office on the fifth floor in Building 24, a drab structure near the center of the MIT.}

\begin{figure}[htb!]
\centering
\includegraphics[scale=.06]{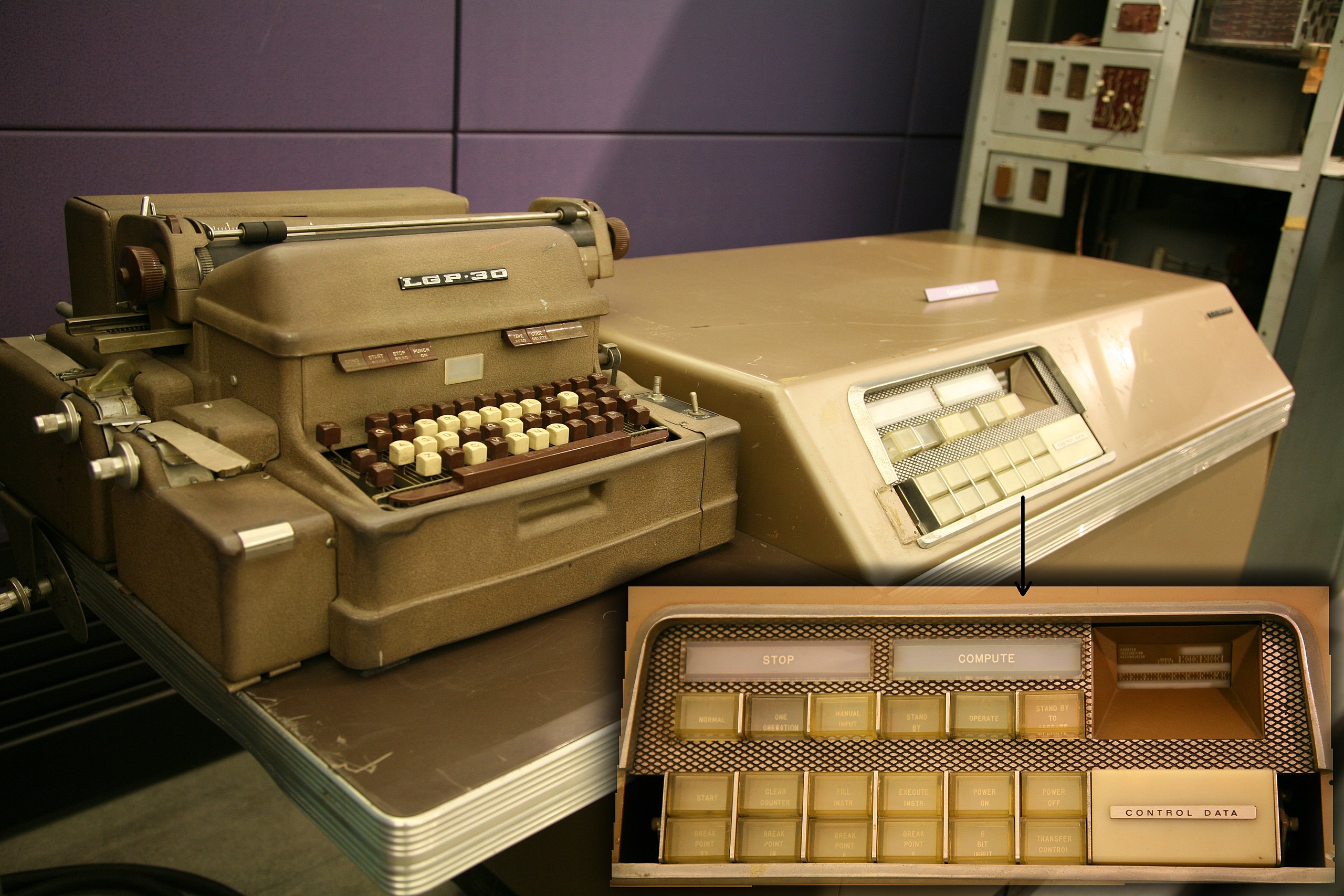}
\legend{\footnotesize{Royal McBee LGP-30}}
%\label{} % a mettre apres caption
\end{figure}

Lorenz published his results in 1963 in a journal of meteorology \cite{lorenz}. For the anecdote, Lorenz's article was sent to Ulam for review. It took 
almost ten years to physicists and mathematicians to realize the importance of this work.
It is only in 1972 that Lorenz presented the ``butterfly effect'' in the 
139th meeting of the American Association for the Advancement of Science, asking: ``Does the flap of a butterfly's wings in Brazil set off a 
tornado in Texas ?''. At this occasion, he presented the surprising picture of the attractor which now bears his name.

Let us point out that it is only in 1998 that Warwick Tucker mathematically proved in his PhD thesis the existence of Lorenz's attractor
\cite{tucker}. His demonstration relies upon a numerical integrator providing a precise control of errors in the approximation of true orbits.

We finish this section by quoting Lorenz \cite{lorenz-1964}:
\begin{displayquote}
\textcolor{unbleu}{
``We thus see that a computing machine may play an important role, in addition to simply grinding out numerical answers. The machine cannot prove 
a theorem, but it can suggest a proposition to be proven. The proposition may then be proven and established as a theorem by analytic means, but 
the very existence of the theorem might not have been suspected without the aid of the machine. Ulam has discussed the general problem of the 
computing machines as a heuristic aid to reasoning, and has presented examples from a number of different branches
of mathematics.''\footnote{Lorenz refers to the book of Ulam \cite{ulam-collection} and to the paper \cite{stein-ulam}.}}
\end{displayquote}

%%%%%%% SECTION
\section{Kruskal, Zabusky, solitons, and the visiometrics}\label{sec:KZ}

\begin{wrapfigure}{r}{0.5\textwidth}
\vspace{-10pt}
  \begin{center}
    \includegraphics[width=0.4\textwidth]{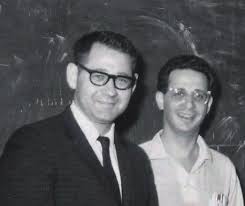}
  \end{center}
  \vspace{-20pt}
  %\caption{\scriptsize{}}
  \legend{{\footnotesize Norman Zabusky and Martin Kruskal.}}
  \vspace{-10pt}
\end{wrapfigure}

At the beginning of the 1960s, in USA, the mathematician and physicist Martin Kruskal (1925-2006) and the physicist Norman Zabusky (1929-2018) started to 
investigate again the FPUT system (see Section \ref{sec:FPUT}). They modified the nonlinear term in the interaction between springs and, together with 
their programmer Gary Deem, they conducted detailed numerical experiments leading to the discovery of a new phenomenon: they observed ``solitary waves''
they called ``solitons'' \cite{ZK}. Remarkably, two solitons can collide with each other and yet preserve their shapes and speeds after collision.
This properties result from a ``symbiotic'' balance of nonlinearity and dispersion.
This numerical experiments were made at the Bell Telephone Laboratory at Whippany on \textsf{IBM 709} and \textsf{7090} computers.

Kruskal and Zabusky realized that a continuous approximation of this system resulted in a partial differential equation which was proposed by 
Diederik Korteweg and his PhD student Gustav de Vries in ...1895!\footnote{This is the equation $\partial_t v+v\partial_{\xi}v+\delta^2\partial^3_{\xi} v=0$.} Their goal was to 
account for the strange waves observed and reported in the 1840s by the Scottish civil engineer John Scott Russel in a canal.
Therefore, Korteweg-de Vries equation was resurrected indirectly by the work of Fermi-Pastal-Ulam-Tsingou, restarted by Kruskal and Zabusky, after 
having been dormant for nearly seventy years. A whole area of physics and mathematics had come into being.

In Kruskal and Zabusky's approach, visualization is of crucial importance. This is of course routine nowadays but it was not the case at that time; 
they had to develop their own visualization tools as well as the tools to interact with the programs. For them, interactive numerical experiments go 
well beyond an assisting tool; the fact that one can modify parameters in real-time and visualize at once the result -- and
then possibly go on -- builds a new intuition and relationship with equations.

We finish this section with two quotes from Zabusky:
\begin{displayquote}
\textcolor{unbleu}{ ``[I] wrote in great detail how Ulam's book\footnote{This is \cite{ulam-collection} (N/A).} had been a great stimulus for me when I discovered it in 1961. At that time I was a lonely runner in an unexplored domain, and Ulam's ideas gave me the fortitude for the long distance I would have to cover. This was the beginning of a new mode of working.''
} (From \cite{Z2005}.)
\end{displayquote}
\begin{displayquote}
\textcolor{unbleu}{``From my observations and the above quotations,\footnote{that are about FPUT model (Section \ref{sec:FPUT}), Feigenbaum work on period-doubling (see Section 
\ref{sec:FCT}), and Lorenz's model (Section \ref{sec:Lorenz} (A/N).} it is clear that interactive and rapid turn-around computing provides an opportunity to concentrate deeply and develop a 
special intuitive ``feel'' for the results. This noninterrupted mode augments the innovative process.''}
(From  \cite{Z1981}.)
\end{displayquote}

%%%%%%% SECTION
\section{H\'enon: from astrophysics to strange attractors}\label{sec:H}

Michel H\'enon (1931-2013) has put numerical experiments at the center of his scientific practice to which he accorded the same status as physics experiments. 
He was interested in astrophysics -- an area where direct experiments are obviously impossible and numerical computations the only way to 
experiment.

\begin{wrapfigure}{r}{0.5\textwidth}
\vspace{-10pt}
  \begin{center}
    \includegraphics[width=0.4\textwidth]{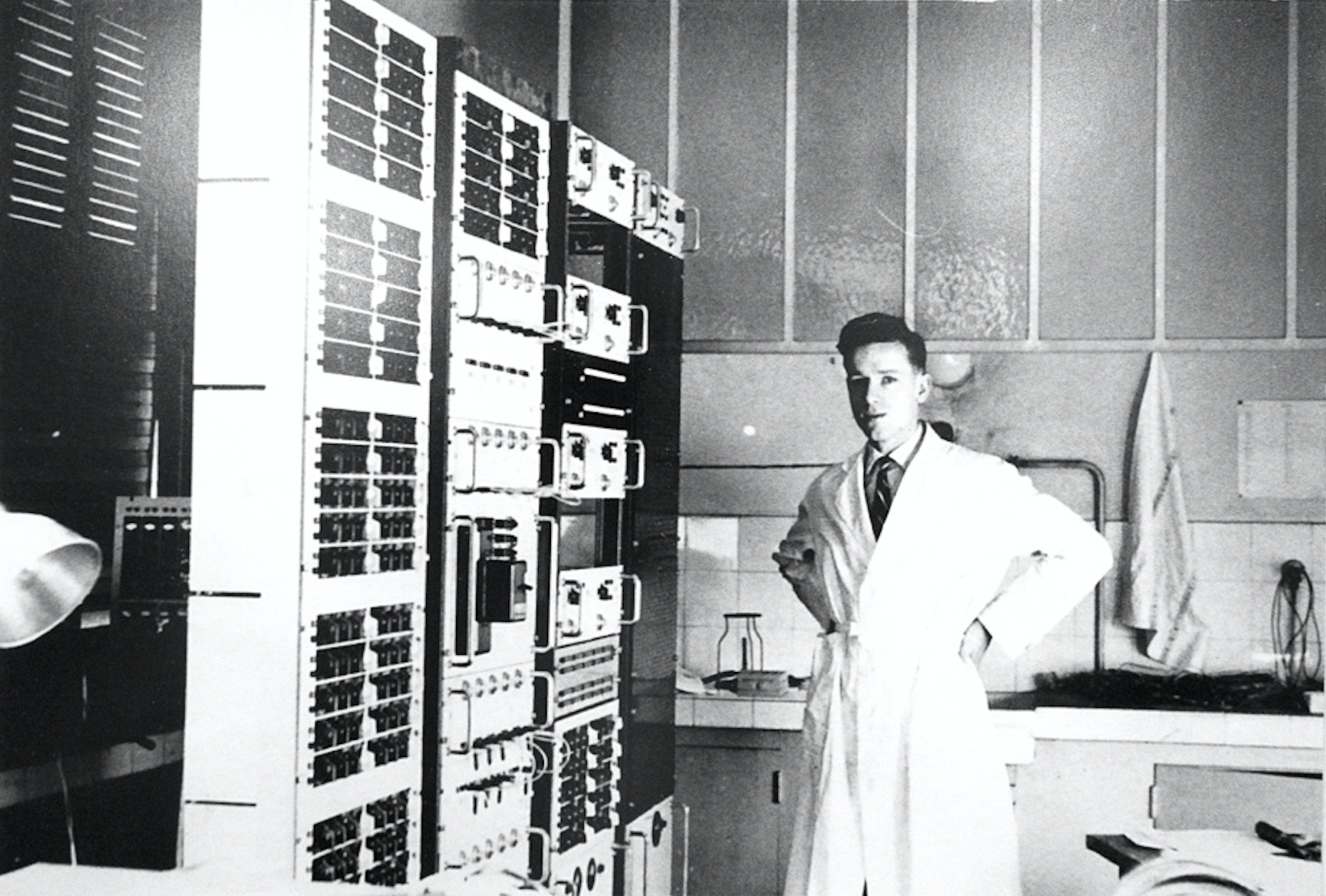}
  \end{center}
  \vspace{-20pt}
  %\caption{\scriptsize{}}
  \legend{{\footnotesize H\'enon and his analogic computer built with Michel Dreux at IAP (circa 1956).}}
  \vspace{-10pt}
\end{wrapfigure}

In the 1950s, H\'enon worked on the construction of analogic computers before making his own one. This was just before the advent and 
democratization of digital computers. At Meudon's observatory, near Paris, he worked on an \textsf{IBM 750} and, later, on an \textsf{IBM 7040} at Nice Observatory. 
He also used the very first programmable pocket calculator, the \textsf{HP-65}.

In the 1960s, H\'enon was interested in various problems arising in astrophysics and, during his 1962 stay at Princeton, he studied the movement of 
a star gravitating in a galaxy with cylindrical symmetry. Several numerical experiments with this model revealed some irregular behaviors. He asked 
Carl Heiles, a graduate student at that time, to redo the program and the experiments by himself on another computer, just as for a physics 
experiments which has to be reproducible.\footnote{They also redid their calculations with different numerical integration methods.}
These experiments resulted in a paper of H\'enon and Heiles in 1964 which revealed a striking mixture of quasi-periodic
and ``chaotic'' behaviors in a model that seemed very simple \cite{HH}, see next figure.

\begin{figure}[htb!]
\centering
\includegraphics[scale=.35]{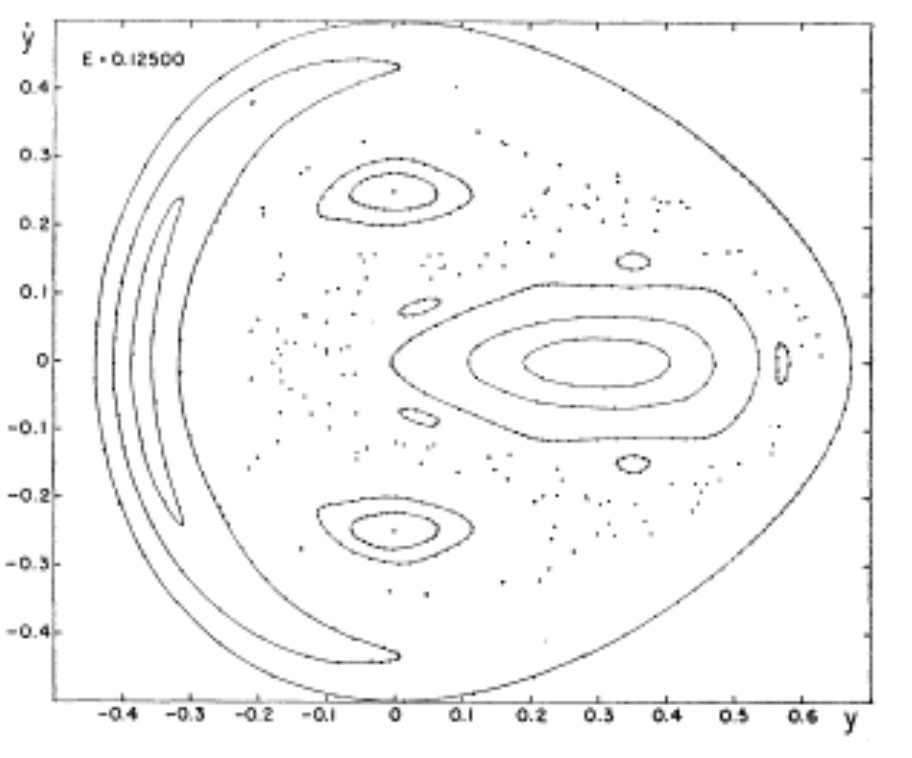}
\includegraphics[scale=.28]{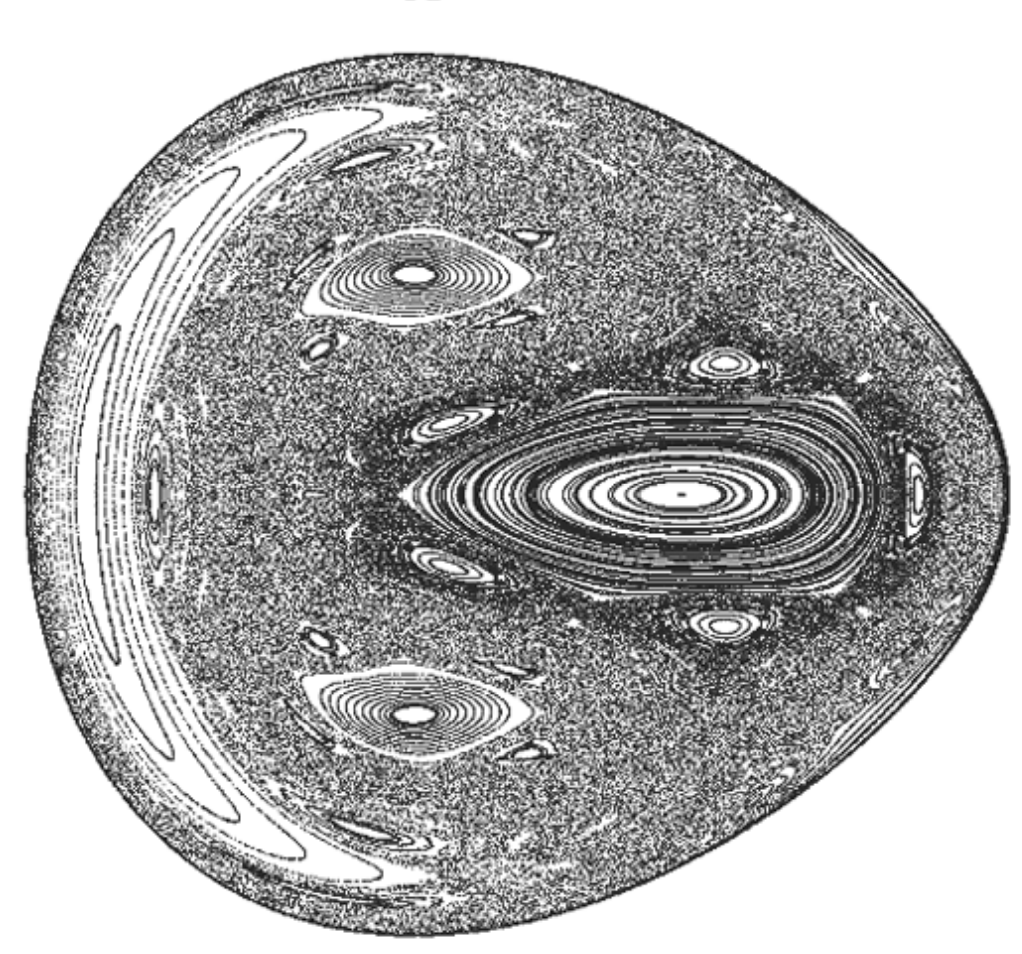}
\legend{\footnotesize{Left: Taken from H\'enon and Heiles paper mentioned above. Right: A version by the authors of this article.}}
%\label{fig:HH}
\end{figure}

H\'enon's approach to astrophysics was to focus on simple mathematical models to shed light on basic phenomena and not to study realistic systems 
which are both intractable and not enlightening. He systematically applied some basic ideas of Henri Poincar\'e and George Birkhoff which lead to 
replace differential equations by iterations of maps using appropriately chosen ``Poincar\'e sections'' to the orbits in phase space.  Although 
seemingly simple, the models so obtained are still very hard to analyze mathematically, but numerical experiments are both easy to make and 
instructive.
A shining illustration of his approach is the so-called ``restricted three body problem'' where he explored systematically the behavior of possible 
orbits \cite{henon-books}. 

\begin{wrapfigure}{r}{0.5\textwidth}
\vspace{-10pt}
  \begin{center}
\includegraphics[width=0.5\textwidth]{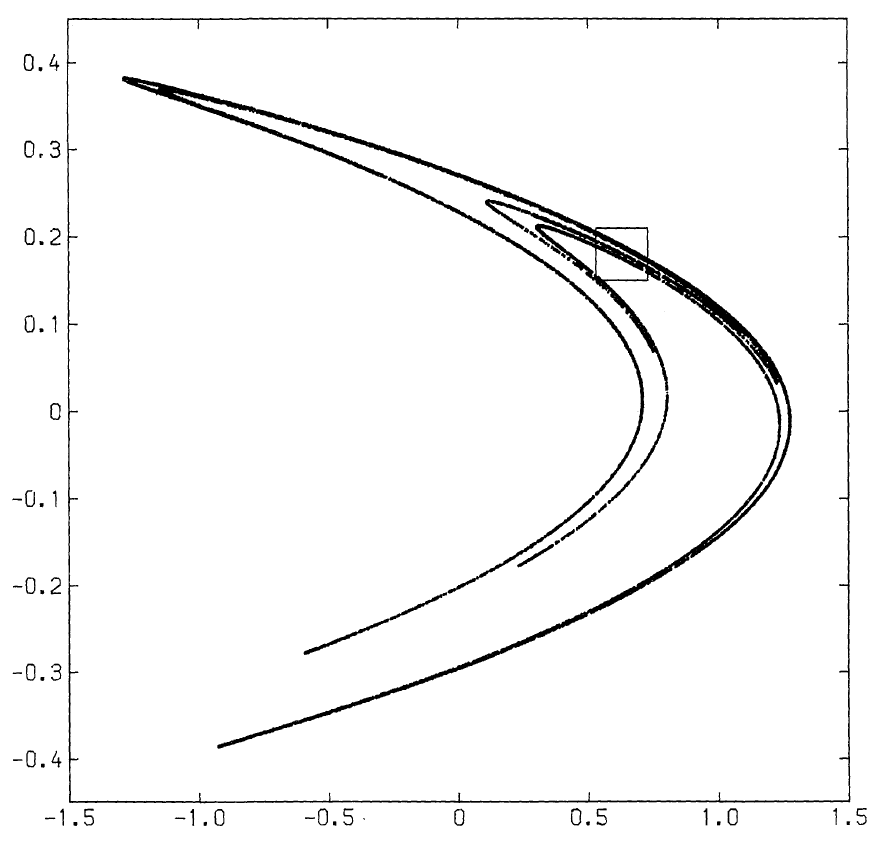}
  \end{center}
  \vspace{-20pt}
  %\caption{\scriptsize{}}
\legend{{\footnotesize From H\'enon's article. \href{https://experiences.mathemarium.fr/Henon-s-attractor.html}{\textsc{click here}} for an online interactive numerical experiment.
}}
  \vspace{-10pt}
\end{wrapfigure}

The work for which H\'enon is undoubtedly famous outside astronomy is the attractor which bears his name.
Tuning the parameters of Lorenz's equations we alluded to above, and using a Poincar\'e section, Yves Pomeau
and Jose-Luis Ibanez observed Smale's horseshoe mechanism. Pomeau gave a talk on this observation at Nice Observatory to which H\'enon assisted. 
This prompted him to propose a very simple model based on a quadratic map in the plane. In this model,
the tuning of a parameter displays Smale's horseshoe mechanism; this is the so-called H\'enon map:
\[
\begin{cases}
x_{n+1}=y_n+1-ax_n^2\\
y_{n+1}=bx_n
\end{cases}
\]
where $a,b$ are positive parameters, and given an initial condition $(x_0,y_0)$. For $a=1.4, b=0.3$, there is a strange attractor displayed in the above figure.

H\'enon made his computations with a \textsf{HP-65} and turned to a \textsf{IBM 7040} computer for more extensive computations.
The fact that it is not a numerical artifact was an open problem till 1991 when the mathematicians Michael Benedicks 
and Lennart Carleson published a (lengthy and complicated) mathematical proof.\footnote{Notice that the ``historical values'' of the parameter are not covered by their 
result.}

%%%%%%% SECTION
\section{The period-doubling cascade of Feigenbaum, Coullet and Tresser}\label{sec:FCT}

In 1978, a young researcher named Pierre Coullet (born in 1949) from the University of Nice and Charles Tresser (born in 1950), a graduate student, got interested in 
the mechanisms leading from laminar to turbulent fluids. Such questions are notoriously hard and they decided to start with very simple mathematical 
models based on iterating a map of the unit interval into itself. As for H\'enon's model, numerical experiments were decisive. The novelty they bring, 
which seems obvious nowadays, is to develop an interactive approach in their experiments; parameters are changed repeatedly in a simple way to 
explore almost in real-time the behavior of the models.

At the same time and independently of Mitchell Feigenbaum (1944-2019), researcher at the Los Alamos Scientific Laboratory, Coullet and Tresser applied 
renormalization group techniques to a class of ``unimodal'' maps (their graph has a single ``bump''). The simplest map in this family is the logistic map.
The value $x_n$ at time $n$ is mapped to a new value $x_{n+1}$ at time $n+1$ as follows: $x_{n+1}=rx_n(1-x_n)$, where the control 
parameter $r$ satisfies $0<r\leq 4$ and $0\leq x_n\leq 1$.
They showed that a transition to chaos arises by a period-doubling bifurcation or ``cascade'' which possesses some universal features \cite{coullet,feigenbaum}. In 
particular, the limiting ratio of each bifurcation interval to the next between every period doubling is a universal constant named the Feigenbaum number.

Using an \textsf{HP-9825} and an XY plotter, Coullet and Tresser could visualize the iterates of the logistic map, change the parameter and see at 
once what  was the resulting effect. They used HPL language (similar to BASIC) so they could change the value of a variable via the keyboard 
without stopping the program. On his side Feigenbaum used the \textsf{HP-65} programmable calculator, the same as Michel H\'enon.

Let us quote Oscar Lanford who proved many of the conjectures of Feigenbaum using computer-assisted proofs (see \cite{lanford}):\footnote{We also refer to the book 
of Collet and Eckmann \cite{collet-eckmann}.}
\begin{displayquote}
\textcolor{unbleu}{
``The methods used to study smooth transformations of intervals are by and large elementary and the theory could have been developed long 
ago \emph{if anyone had suspected that there was anything worth studying.} In actual fact, the main phenomena were discovered through
numerical experimentation and the theory has been developed to account for the observations. In this respect, computers have played a crucial
role in its development.'' \cite{lanford-bourbaki}
}
\end{displayquote}

Is it possible to observe transition to chaos by a period-doubling cascade in a physics experiment or is it purely mathematical phenomenon? 
The answer is positive. In 1979, Libchaber, Laroche and Fauve \cite{libchaber} observed this scenario in a beautiful convection experiment using liquid helium and later 
mercury. Tuning the temperature difference between the bottom and the top of the box containing the fluid, he observed the successive bifurcations with an 
amazing precision.

\begin{figure}[htb!]
\centering
\includegraphics[scale=.18]{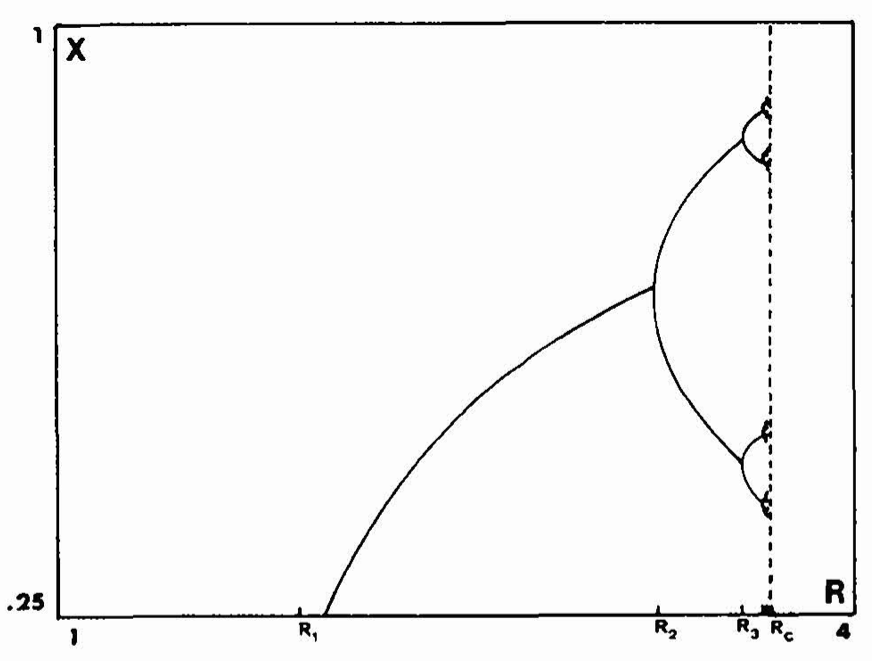}
\includegraphics[scale=.052]{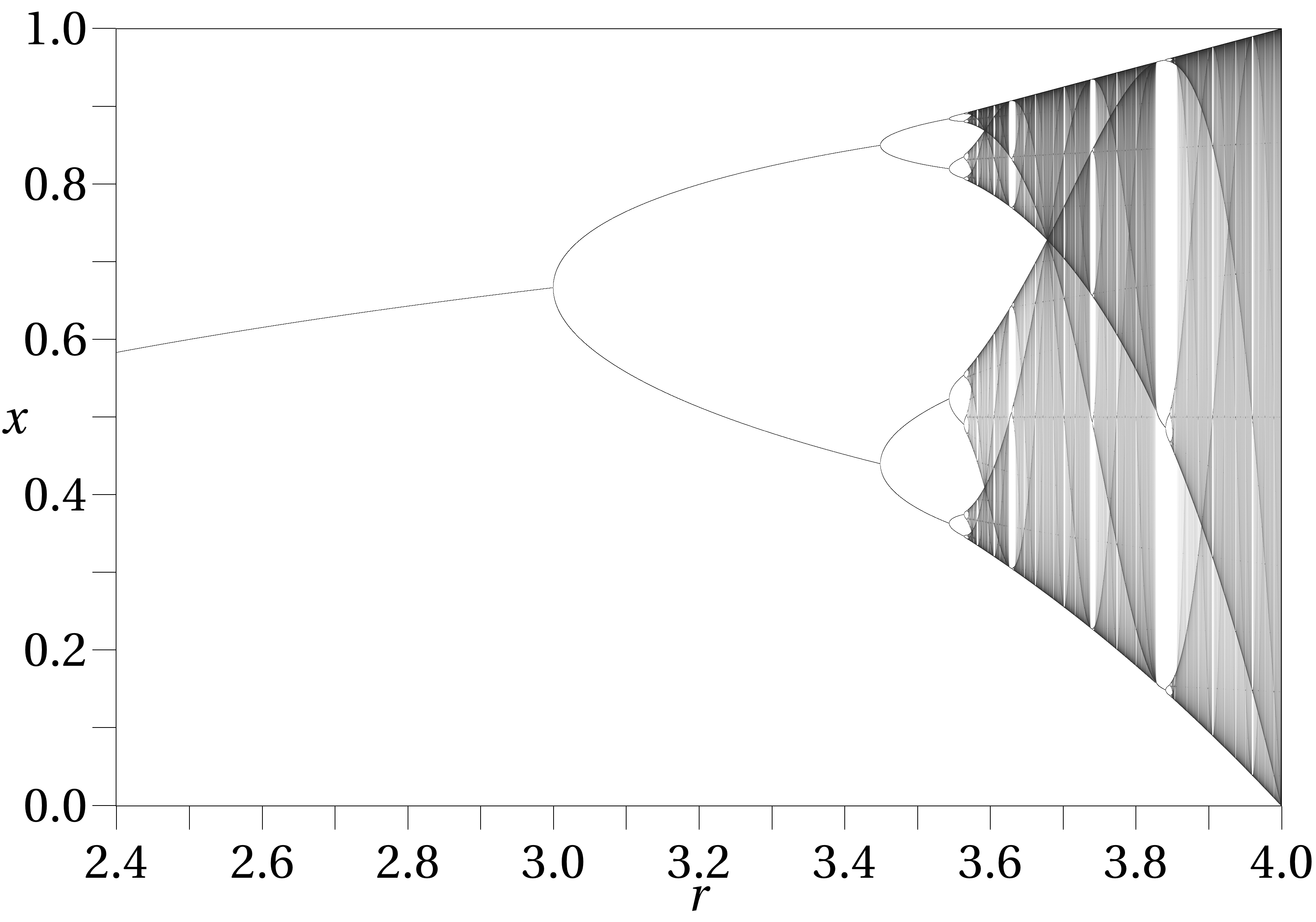}
\legend{\footnotesize{Bifurcation diagram for the logistic map. The attractor for any value of the parameter $r$ is
shown on the vertical line at that $r$. On the left: from Coullet-Tresser article (in which they use $X$ and $R$ instead of $x$ and $r$). On the right: from \href{https://en.wikipedia.org/wiki/Logistic_map}{Wikipedia}.\newline \href{https://experiences.mathemarium.fr/Logistic-map.html}{\textsc{click here}} for an online interactive numerical experiment with which you can explore the behaviour of the orbits and observe the appearance of chaos, and \href{https://experiences.mathemarium.fr/Bifurcation-diagram-of-the.html}{\textsc{click here}} for the bifurcation diagram in which you can zoom.}}
%\label{} % a mettre apres caption
\end{figure}

%%%%%%% SECTION
\section{Iterating complex polynomials:  Hubbard and Mandelbrot}\label{sec:HM}

We end this article with complex dynamics, a field which was pioneered by the French mathematicians Pierre Fatou and
Gaston Julia at the beginning of the 20th century. The fact that it became dormant for about sixty years is not by accident: without 
the visualization of Julia sets provided by computers, it is obvious that nobody had  only the haziest notion of what these sets might look like when 
plotted on the plane.
  
\begin{figure}
\centering
\includegraphics[width=0.3\textwidth]{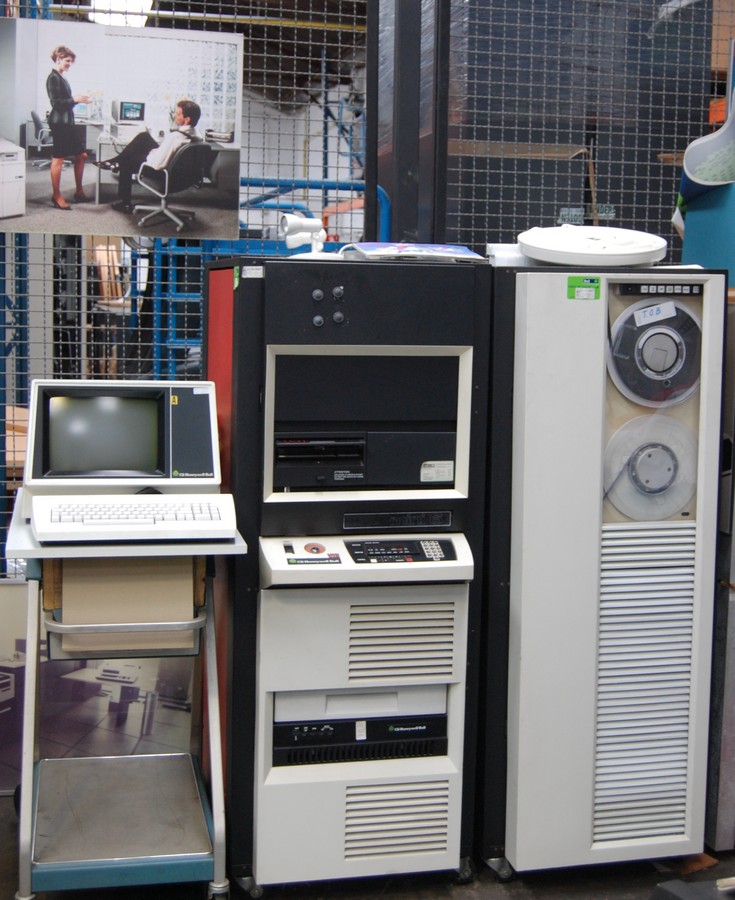}\hspace{1cm}
\includegraphics[width=0.35\textwidth]{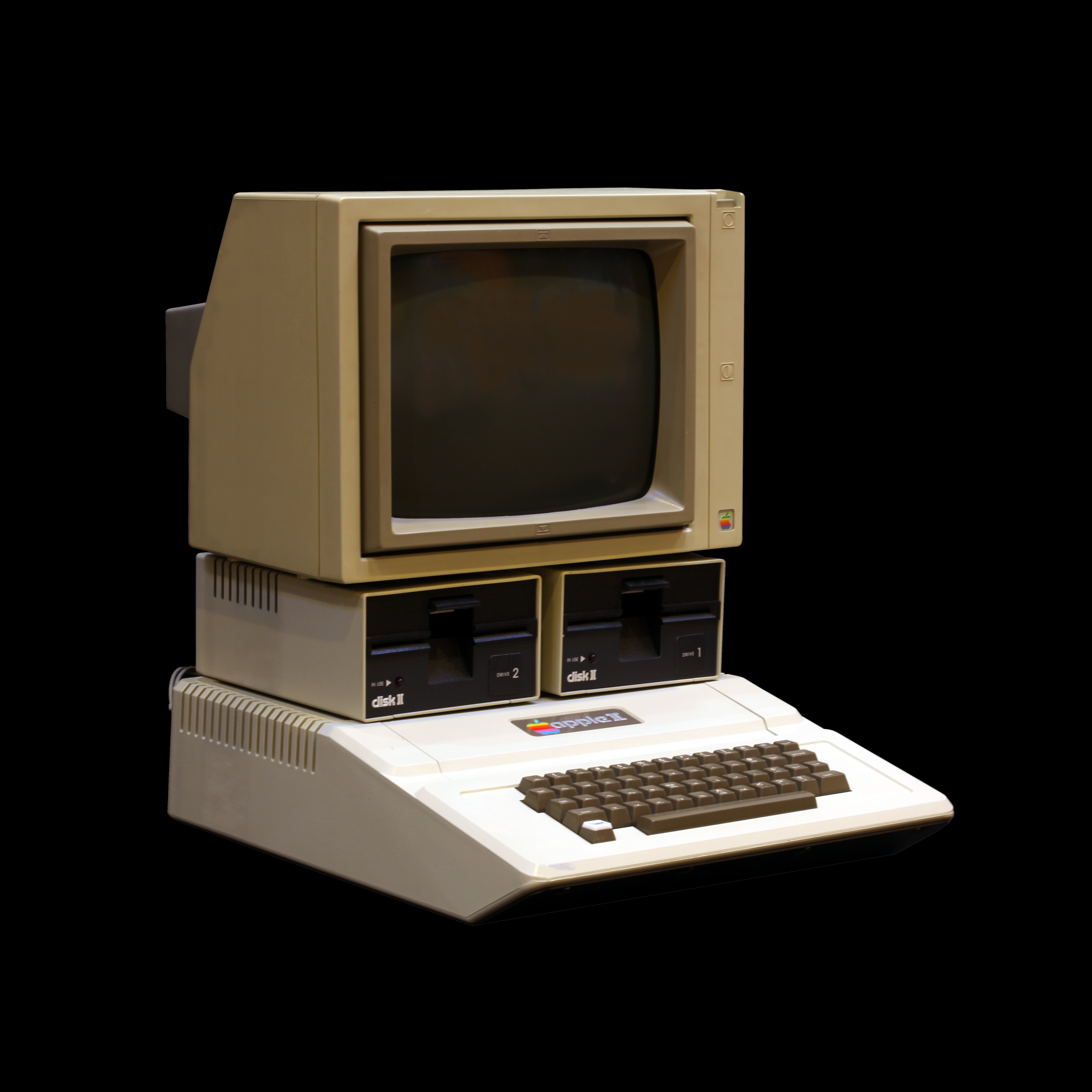}
\legend{On the left: \textsf{Mini 6} (Honeywell-Bull, 1978). On the right:  \textsf{Apple II (1977)}}
\label{fig-mini6}
\end{figure}

In the preface of the book ``The Mandelbrot set, theme and variations'' \cite{hubbard}, the mathematician John Hubbard (born in 1945) explained
how teaching undergraduate students pushed him to use numerical experiments:
\begin{displayquote}
\textcolor{unbleu}{
During the academic year 1976-77, I was teaching DEUG (first and second year calculus) at the University of Paris at Orsay. At the time it was
clear that willy-nilly, applied mathematics would never be the same again after computers. So I tried to introduce some computing into the
curriculum. [...] Casting around for a topic sufficiently simple to fit into the 100 program steps and eight memories of these primitive machines, 
but still sufficiently rich to interest the students, I chose Newton's method\footnote{Applied to a complex polynomial, like $z^3-1$ (A/N).} for 
solving equations (among several others). [...] But when a student asked me how to choose an initial guess, I couldn't answer. It took me some 
time to discover that no one knew, and even longer to understand that the question really meant: what do the basins of the roots look like?''
}
\end{displayquote}
%at Orsay University, in Paris, during the year 1976-77, pushed him to use numerical
%experiments. Casting around for a sufficiently simple topic fitting the first and second year calculus course, but implementable
%of the primitive machines available, he chose Newton's method for solving equations. Working in complex analysis, he wanted to apply it to a 
%complex polynomial, $z^3-1$ for instance. He realized that in fact no one knew what the basins of the roots looked like!
He continues:
\begin{displayquote}
\textcolor{unbleu}{
``As I discovered later, I was far from the first person to wonder about this. Cayley \cite{cayley} had asked about it explicitly in the 1880's, and Fatou and Julia had explored some cases around 1920. But now we could effectively answer the question: computers could draw the basins. And they did: the math department at Orsay owned a rather unpleasant computer called a \textsf{mini-6} (see Fig.\ref{fig-mini6}), which spent much of the spring of 1977 making such computations, and printing the results on a character printer, with X's, 0's and 1's to designate points of different basins. Michel Fiollet wrote the programs, and I am extremely grateful to him, as 1 could never have mastered that machine myself.''
}
\end{displayquote}

\begin{wrapfigure}{r}{0.48\textwidth}
\vspace{-20pt}
  \begin{center}
    \includegraphics[width=0.25\textwidth]{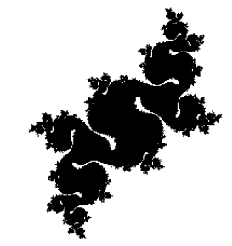}
    \includegraphics[width=0.22\textwidth]{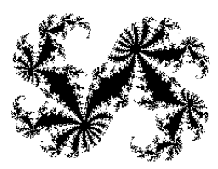}\\
    \includegraphics[width=0.18\textwidth]{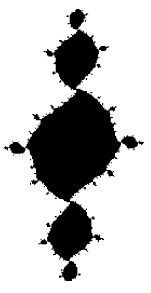}
    \includegraphics[width=0.205\textwidth]{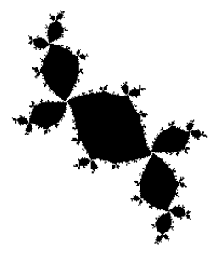}
  \end{center}
  \vspace{-15pt}
  %\caption{\scriptsize{}}
  \legend{{\footnotesize A few Julia sets. \href{https://experiences.mathemarium.fr/Julia-sets.html}{{\textsc{click here}}} for an online interactive
  digital experiment where you can tune $c$.}}
  \vspace{-30pt}
\end{wrapfigure}

With the help of Michel Fiollet, he made some numerical experiments to explore these basins. Stimulated by the mathematician
Dennis Sullivan who was in residence at IHES (Bures-sur-Yvette\footnote{Bures-sur-Yvette is a small town near Paris.}), he turned to plot various Julia sets: take two complex
numbers $z_0$ and $c$ and define the sequence $(z_n)$ by recurrence by setting $z_{n+1}=z_n^2+c$.
For a fixed value of $c$, the Julia set is the frontier of the set of initial values $z_0$ such that the sequence remains bounded.

It seems that Hubbard showed his pictures to Beno\^{\i}t Mandelbrot in 1977, during a stay in USA.
Mandelbrot  told him having often thought about this kind of sets although he had never
made pictures of them. Hubbard mentions in the above-mentioned text that the arrival of the \textsf{Apple II} was decisive in 
making better pictures in an easier way. 

Although Mandelbrot (1924-2010) worked at IBM and had access to the best computers available at that time, it is during 
a stay at Harvard that he obtained for the very first time, in March 1980, a rough image of the set
which now bears his name -- the Mandelbrot set. He used a \textsf{Vax} computer. Peter Moldave, a teaching assistant, 
volunteered his services as a programmer. Let us quote Mandelbrot:
\begin{displayquote}
\textcolor{unbleu}{
[...] I discovered the Mandelbrot set in 1980 while at Harvard, at a time when the computer facilities there were among the most miserable in academia. The basement of the
Science Center housed its first D.E.C. Vax 50 (not yet `broken in'), pictures were viewed on a Tektronix cathode-ray (worn out and very faint), and hard copies were printed on
a ill-adjusted Versatec device. We could only work at night when we had only one competitor for the machine, chemist Martin Karplus. 
\cite[p. 21]{mandelbrot-book}
}
\end{displayquote}

To obtain this set, one plots the set of all the values of $c$ such that the above-defined sequence
remains bounded, starting from $z_0=0$ each time.  Mandelbrot made much better pictures of it at IBM
facilities still with the aid of Moldave and published an article in 1980 reporting this discovery. It is interesting to quote Mandelbrot again
\cite{mandelbrot-book}:
\begin{displayquote}
\textcolor{unbleu}{
``In any event, IBM was absolutely not graphics-oriented. I had no roomful of up-to-date custom equipment, only good friends who built a custom 
contraption that established the state of the art in image rendering.''
}
\end{displayquote}
The mathematical study really started in 1984 with the work of Adrien Douady and Hubbard who
established the fundamental properties of this extraordinary set and named it after Mandelbrot \cite{douady-hubbard-1,douady-hubbard-2}.
Hubbard made many numerical experiments to guide their intuition.
In 1985, the mathematicians Heinz-Otto Peitgen and Peter Richter popularized the Mandelbrot set by making striking colorful pictures of it \cite{beauty}.

\begin{figure}[htb!]
\centering
\includegraphics[scale=.12]{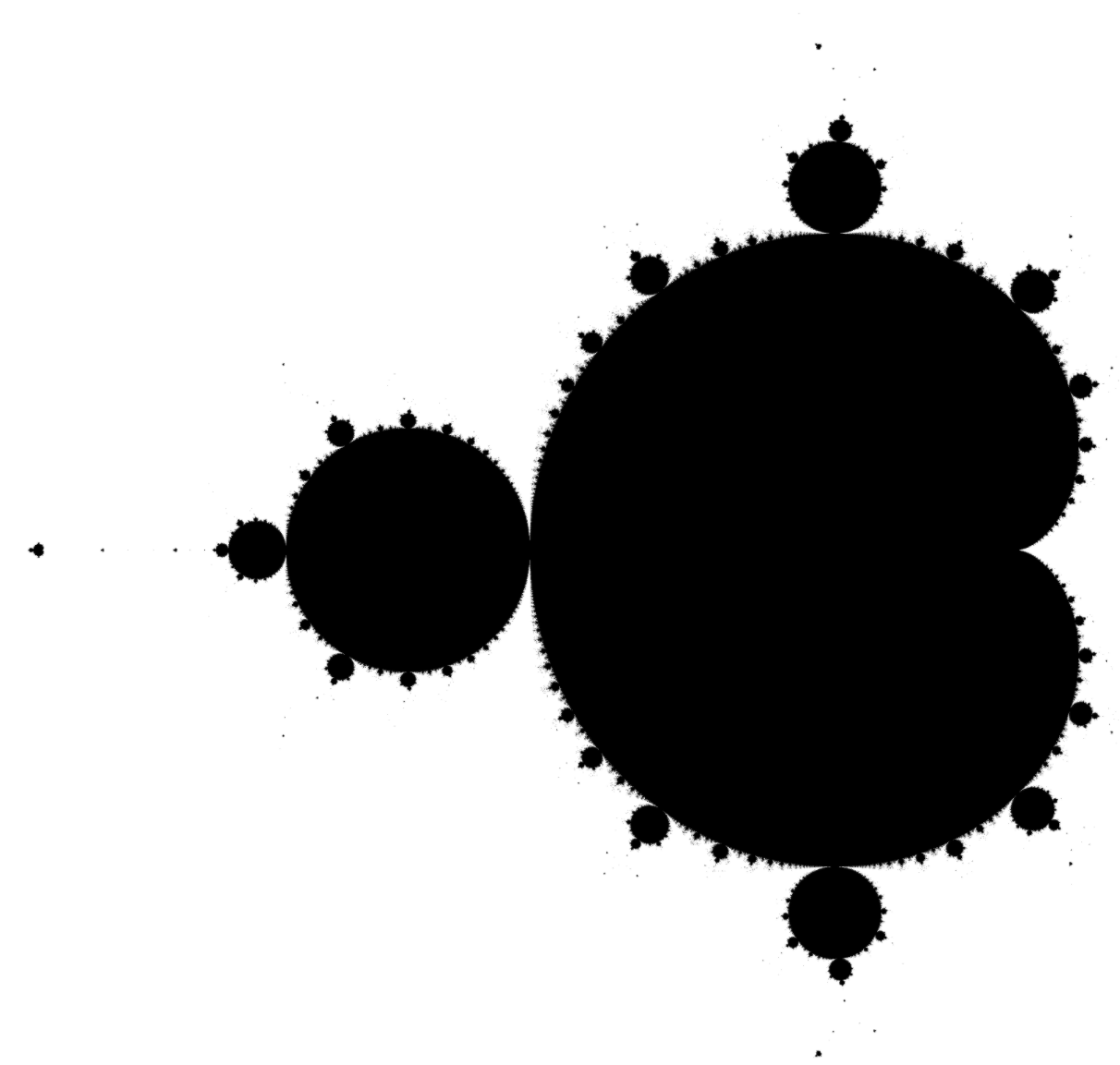}\hspace{.5cm}
\includegraphics[scale=.15]{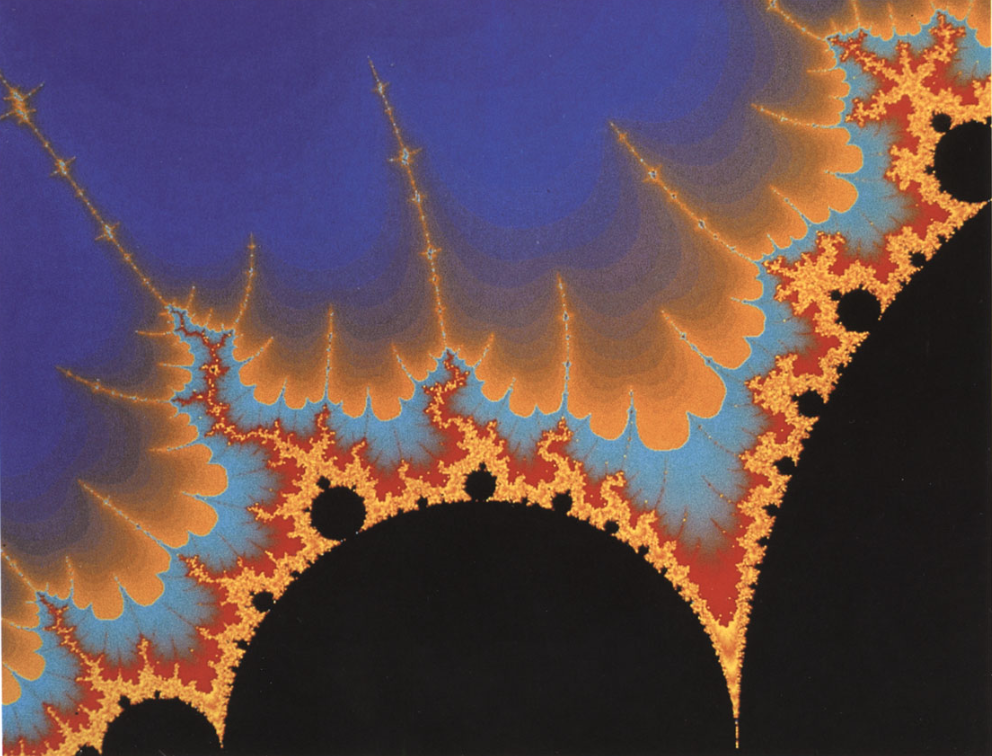}
\legend{\footnotesize{Left:  The Mandelbrot set. Right: Part of the Mandelbrot set, taken from \cite{beauty}.}}
%\label{fig:HH}
\end{figure}

By the way, we recommend the following video-lecture of Hubbard:
\href{https://www.youtube.com/watch?v=NKZIm70jYjQ}{The Beauty and Complexity of the Mandelbrot Set},  dating from 1989.\footnote{This is a 
videotape which can be found in dvd format on the \href{https://bookstore.ams.org/dvd-27}{AMS book store}. We found it on Youtube.}

\newpage 

%%%%%%% SECTION

%\addcontentsline{toc}{section}{Appendix: More pictures}

\section{More pictures}\label{sec:more-pics}

We have collected a small selection of pictures illustrating the preceding material. Needless to say that we have omitted a great deal of pictures that we would have liked to include.

%\begin{figure}
\begin{figure}[htb!]
\centering
\includegraphics[scale=.6]{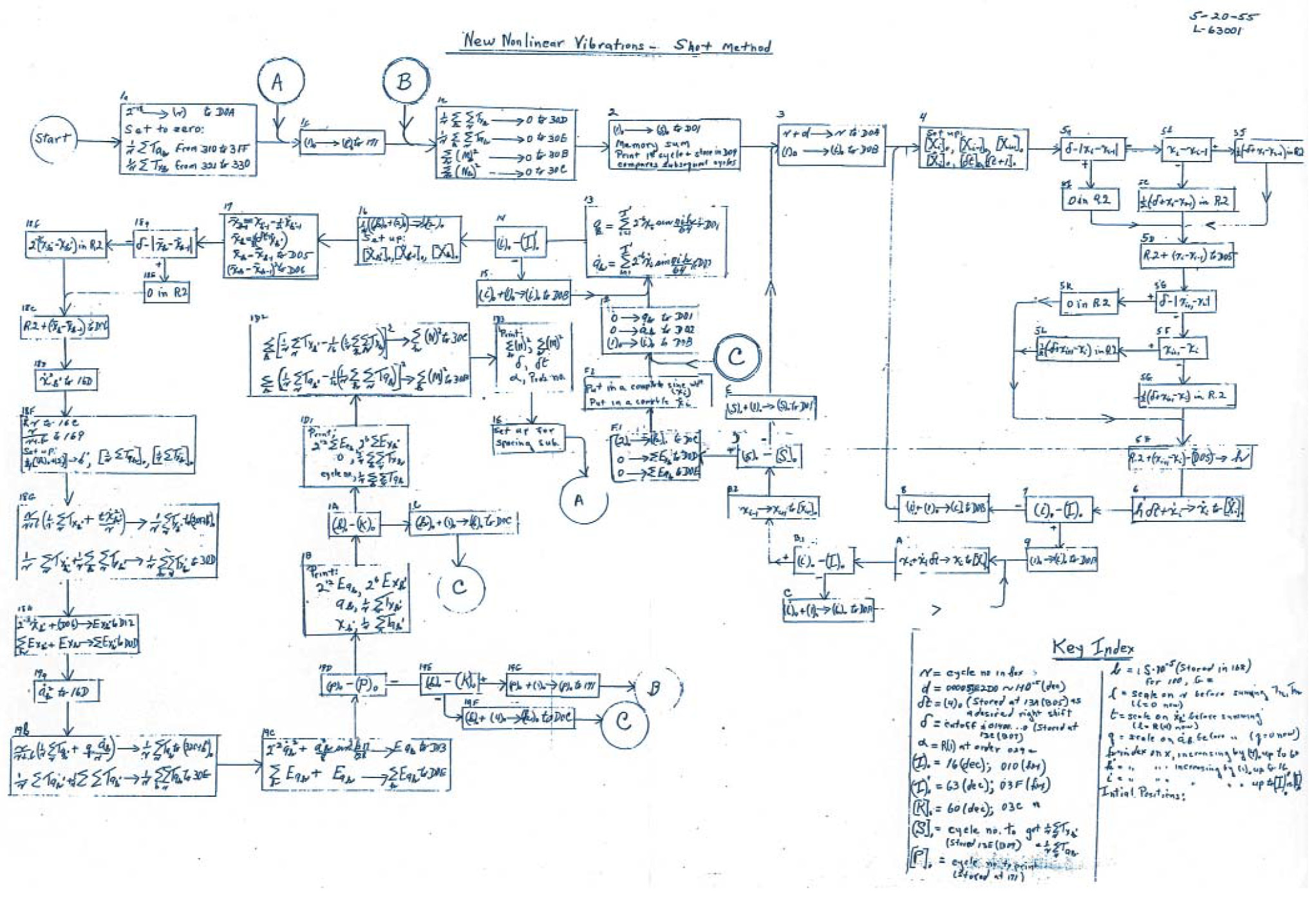}
\legend{Reproduction of the algorithm used by Mary Tsingou to similate the FPUT model (see Section \ref{sec:FPUT}).}
%\label{} % a mettre apres caption
\end{figure}

\begin{figure}[htb!]
\centering
\includegraphics[scale=.4]{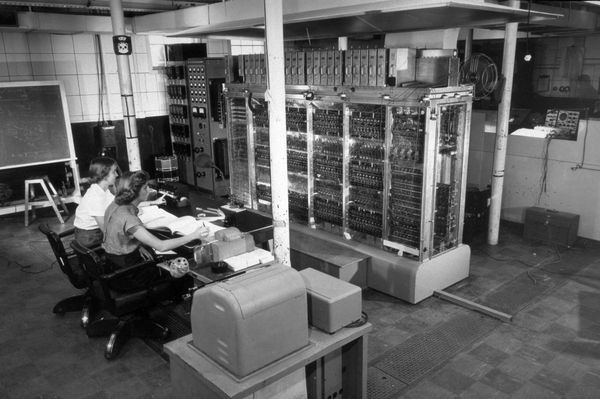}
\legend{The MANIAC I (Mathematical Analyzer Numerical Integrator and Automatic Computer Model I) was an early computer built under the direction of Nicholas Metropolis at the Los Alamos Scientific Laboratory. It was based on the von Neumann architecture of the IAS, developed by John von Neumann. MANIAC I. In 1953, the MANIAC obtained the first equation of state calculated by modified Monte Carlo integration over configuration space. Source: \href{https://en.wikipedia.org/wiki/MANIAC_I}{Wikipedia}.}
%\label{} % a mettre apres caption
\end{figure}

\begin{figure}
\centering
\includegraphics[scale=.15]{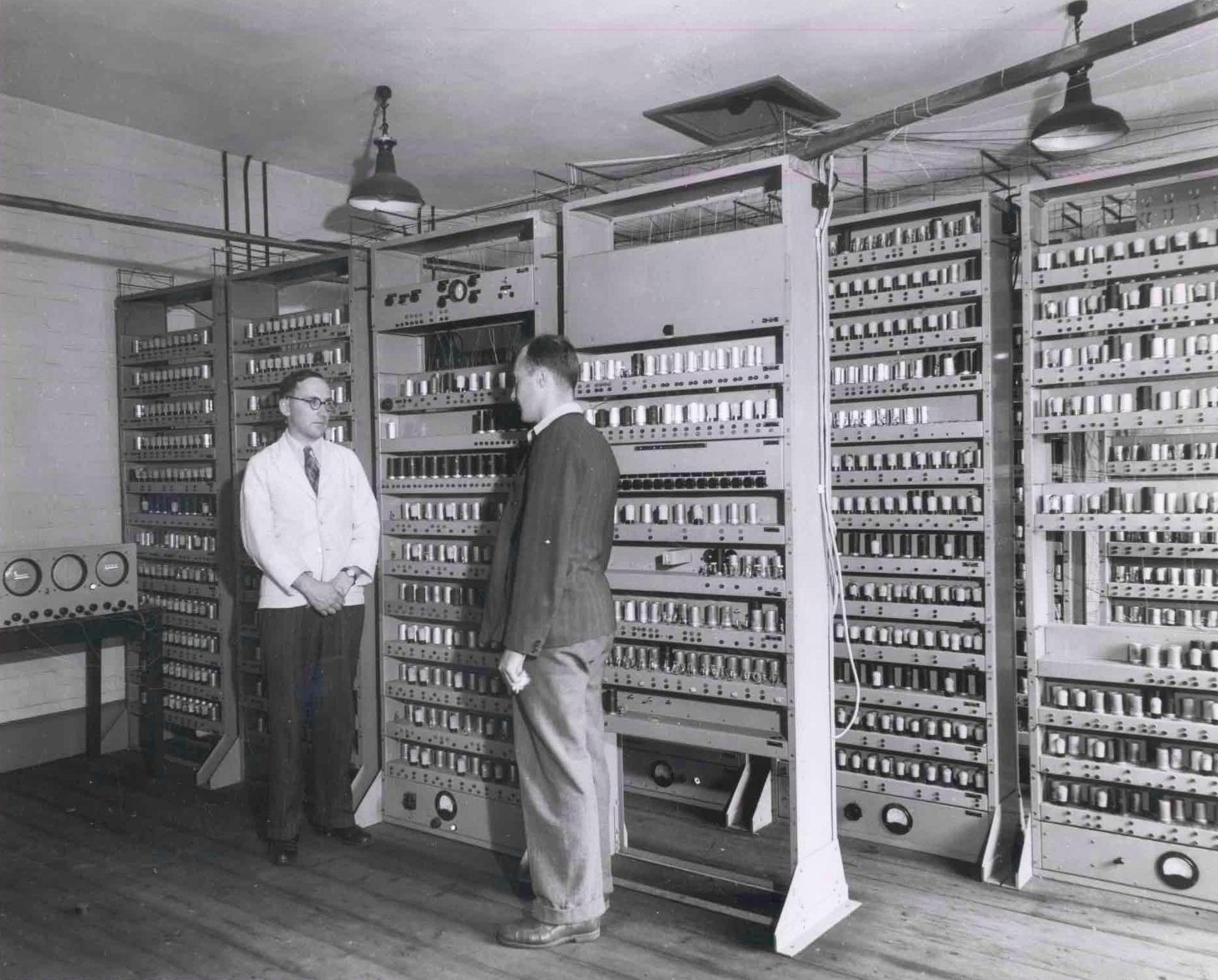}
\legend{The Electronic delay storage automatic calculator (EDSAC) was an early British computer Inspired by von Neumann's seminal First Draft of a Report on the EDVAC, the machine was constructed by Maurice Wilkes and his team at the University of Cambridge Mathematical Laboratory. It was released in 1949. Source: \href{https://en.wikipedia.org/wiki/EDSAC}{Wikipedia}.}
%\label{} % a mettre apres caption
\end{figure}

\begin{figure}
\centering
\includegraphics[scale=.22]{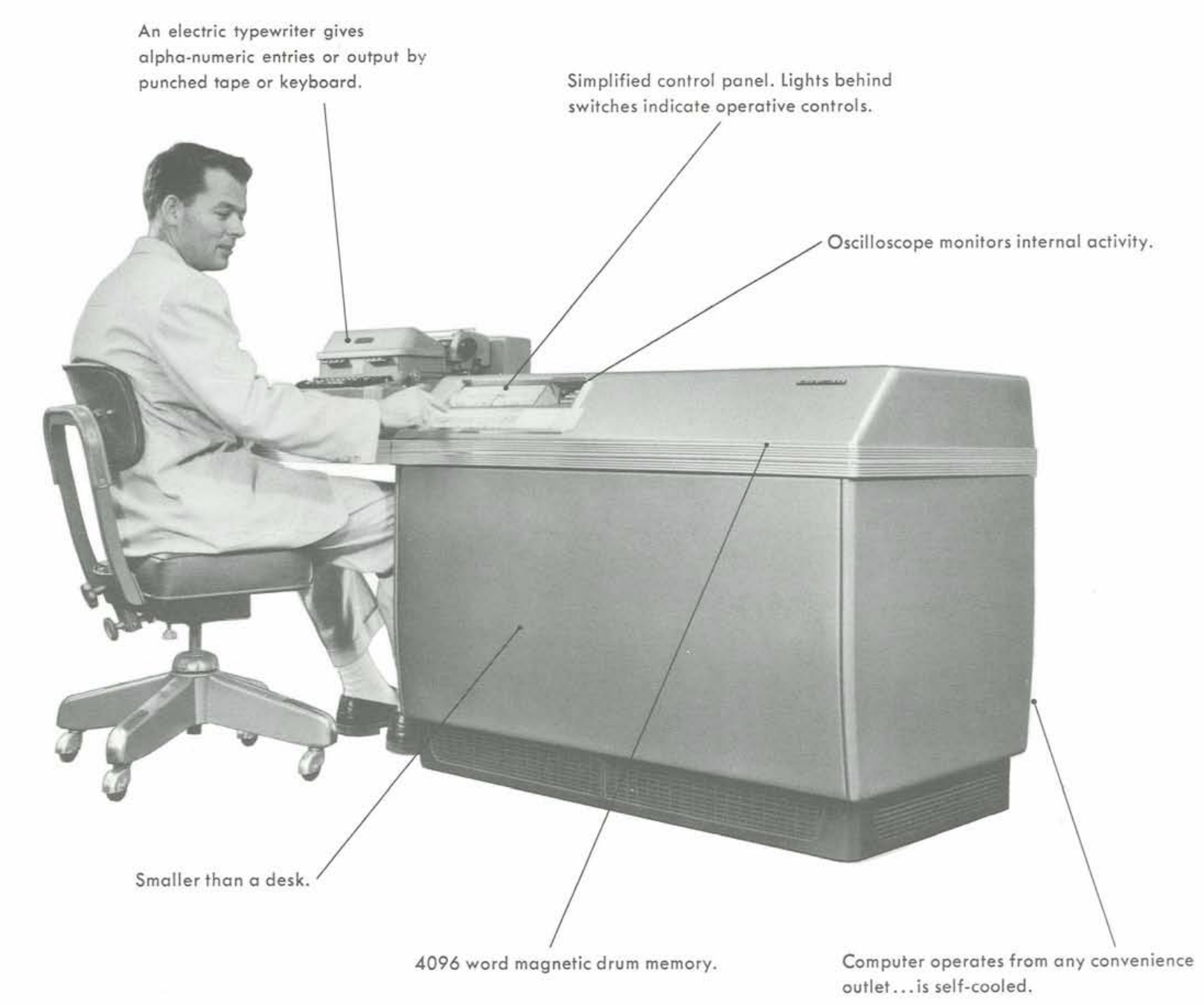}\hspace{0.2cm}
\includegraphics[scale=.22]{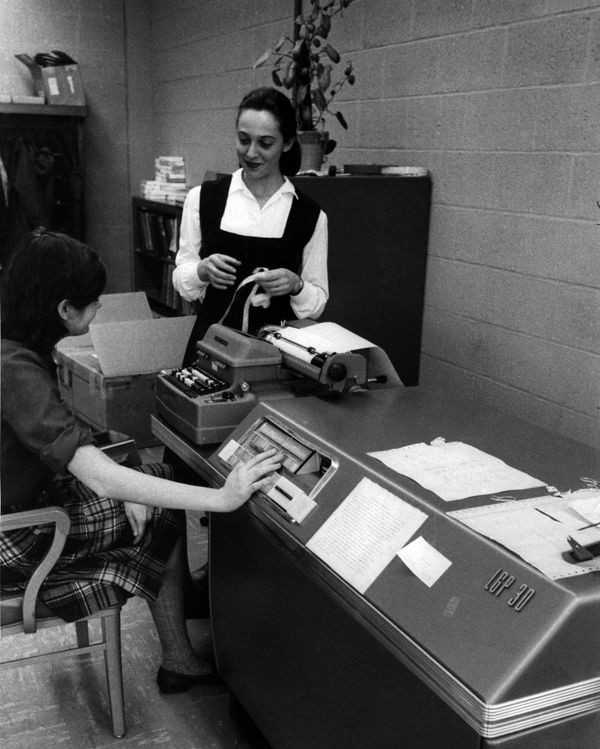}
\legend{\href{https://en.wikipedia.org/wiki/LGP-30}{Royal McBee LGP-30} (1956)}
%\label{} % a mettre apres caption
\end{figure}

\begin{figure}
\centering
\includegraphics[width=0.6\textwidth]{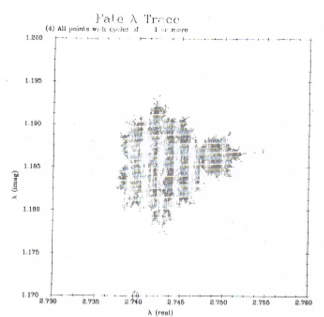}
\legend{\textsf{Mandelbrot set} obtained by Mandelbrot and Moldave in 1 March 1980.}
\end{figure}

\begin{figure}
\centering
\includegraphics[width=0.4\textwidth]{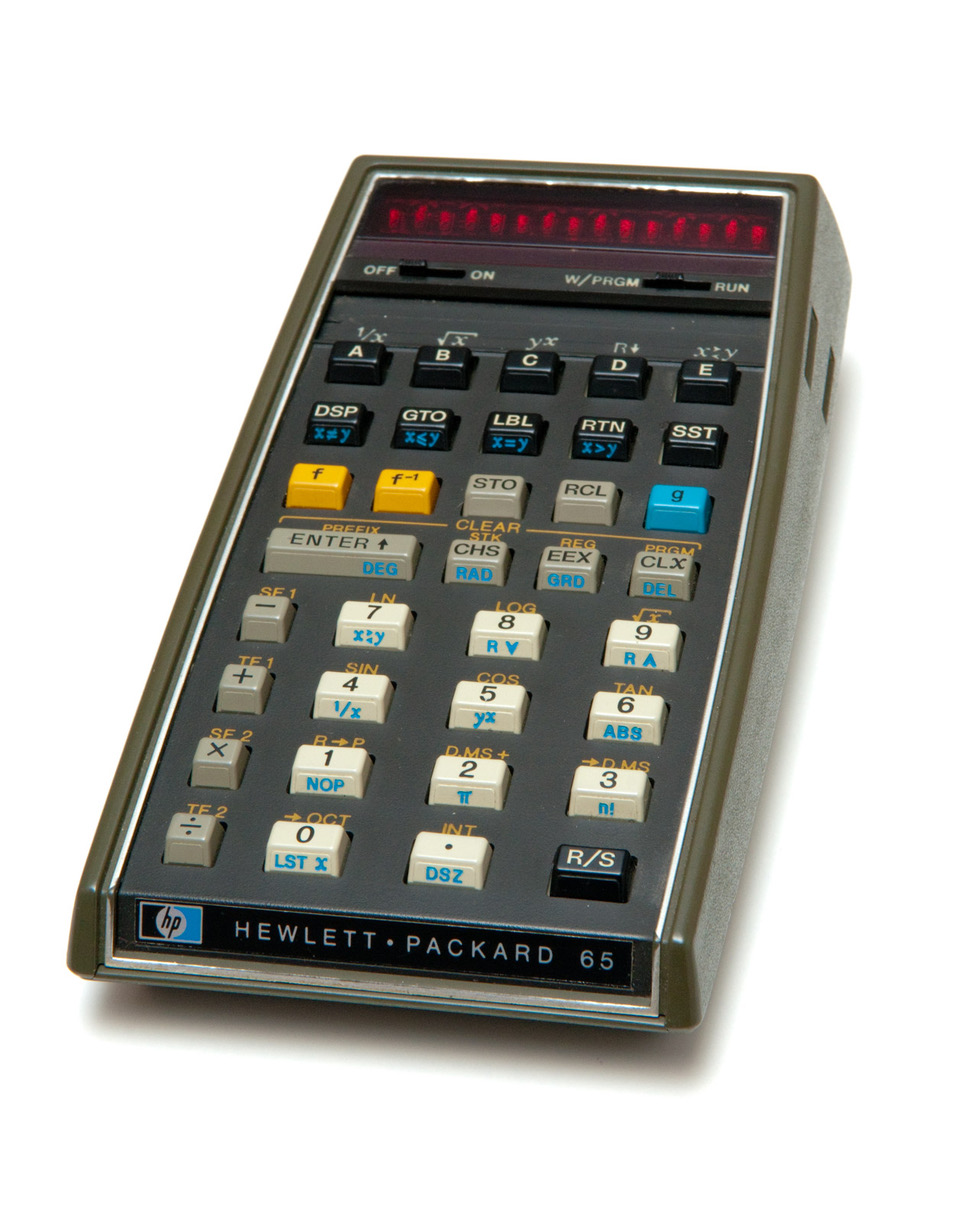}
%\caption{\scriptsize{\textsf{HP-65}}}
\legend{\textsf{HP-65} (released in 1974)}
\end{figure}

\begin{figure}
\centering
\includegraphics[width=0.5\textwidth]{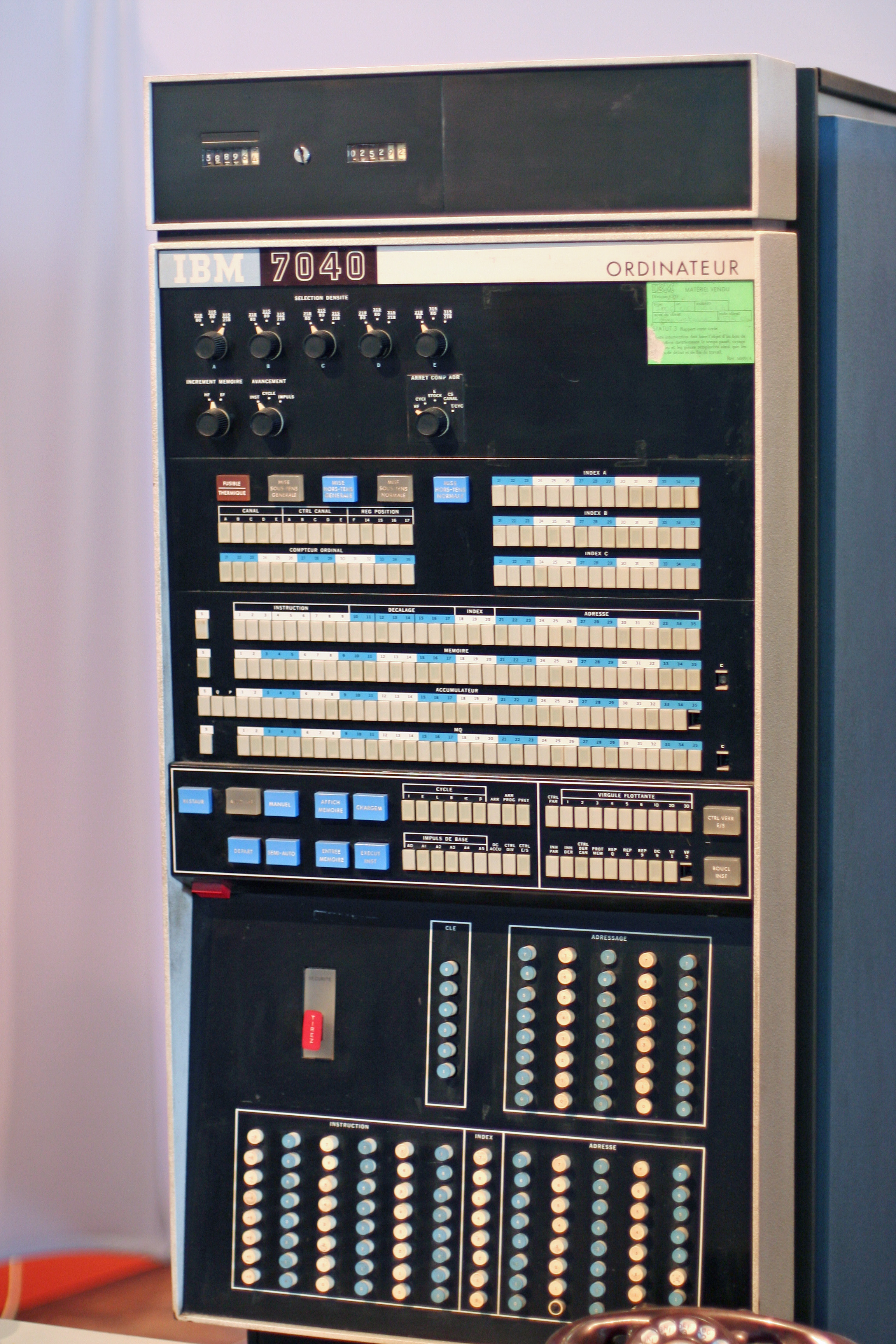}
\legend{Front panel of an \textsf{IBM 7040}}
\end{figure}

\begin{figure}
\centering
\includegraphics[scale=.3]{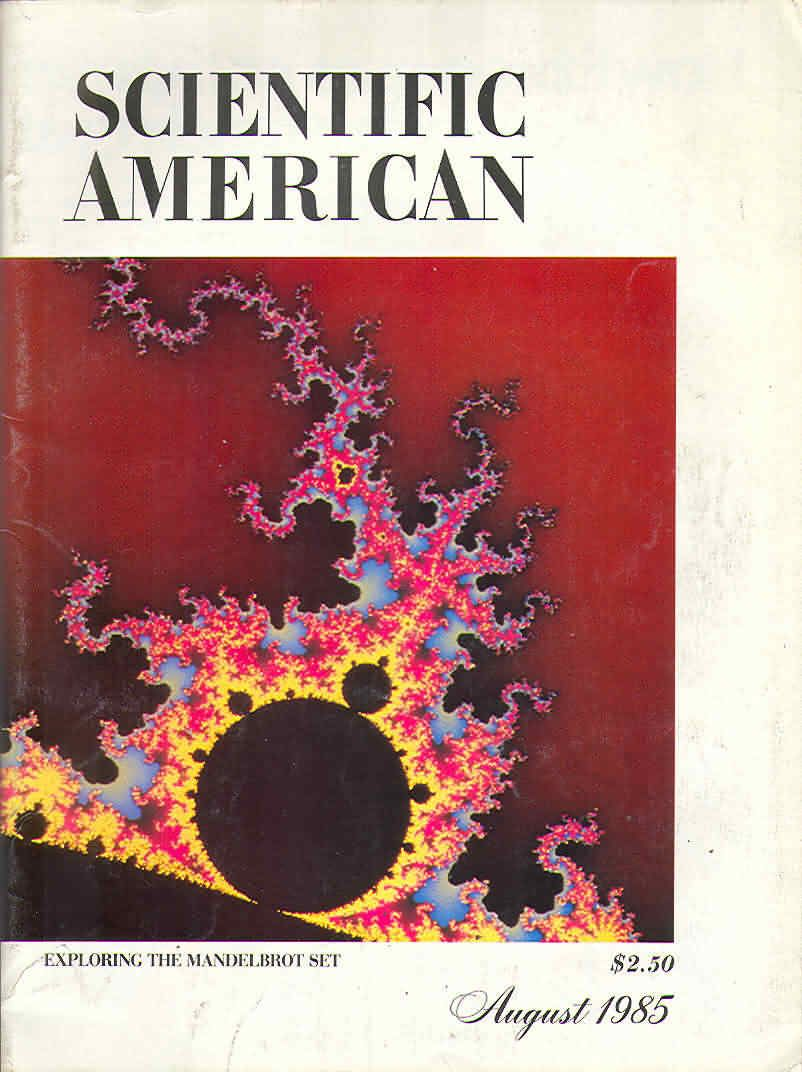}
\legend{}
%\label{} % a mettre apres caption
\end{figure}

\begin{figure}
\centering
\includegraphics[scale=.1]{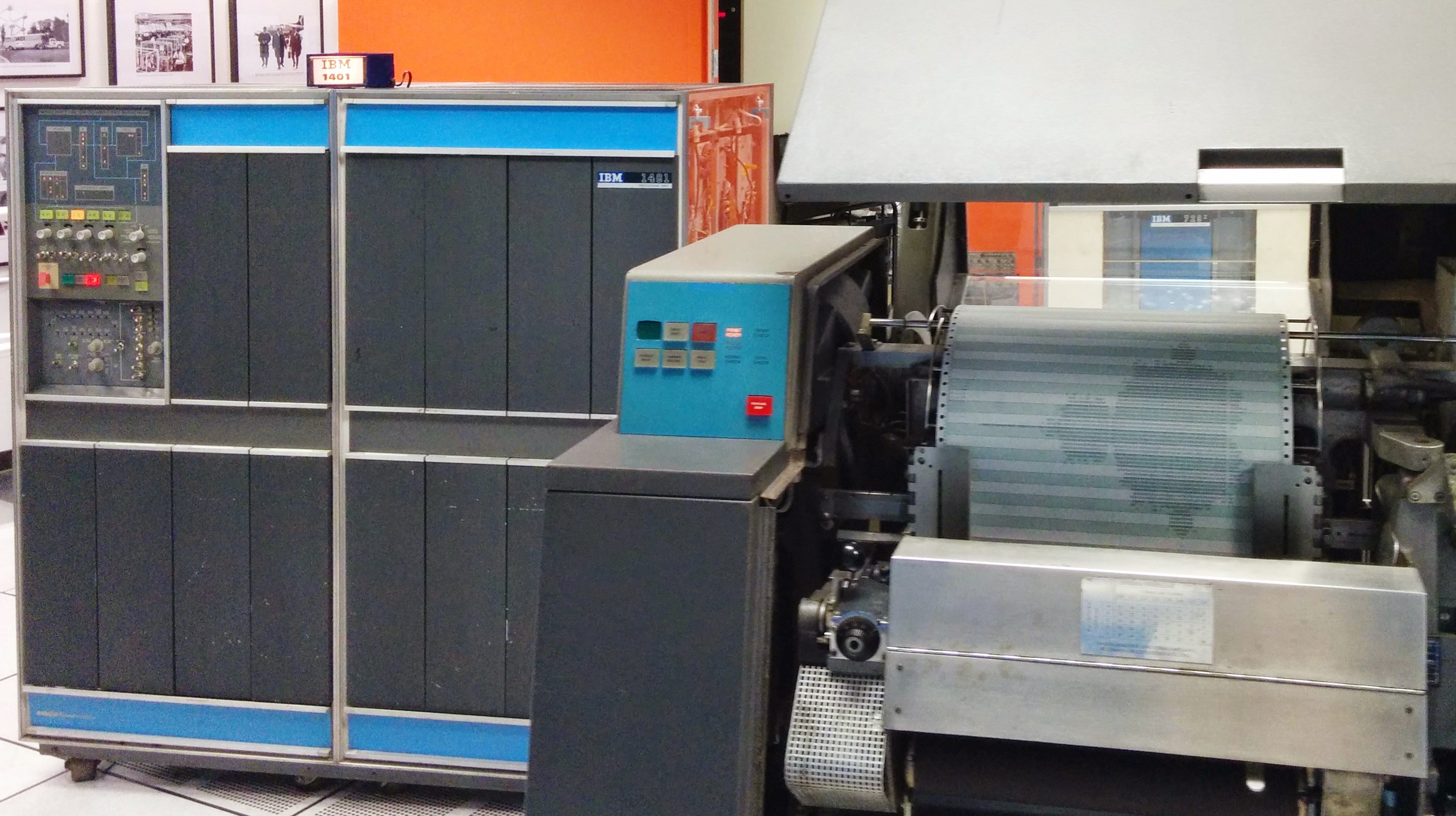}
\legend{The IBM 1401 mainframe computer (left) at the Computer History Museum printing the Mandelbrot fractal on the 1403 printer (right).
Photography taken from \href{http://www.righto.com/2015/03/12-minute-mandelbrot-fractals-on-50.html}{K. Shirriff's blog}.
}
%\label{} % a mettre apres caption
\end{figure}

\newpage

%%%%%%%%%%%%%%%%%%%%%%%%%%%%%%%%%%%%%%%%%%%%%%%%%%%%%%%%%%%%%%%%%%%%%%
%%%%%%%%%%%%%%%%%%%% BIBLIO
%%%%%%%%%%%%%%%%%%%%%%%%%%%%%%%%%%%%%%%%%%%%%%%%%%%%%%%%%%%%%%%%%%%%%%

\end{document}